\documentclass{article}

\usepackage{PRIMEarxiv}

\usepackage[utf8]{inputenc} % allow utf-8 input
\usepackage[T1]{fontenc}    % use 8-bit T1 fonts
\usepackage{hyperref}       % hyperlinks
\usepackage{url}            % simple URL typesetting
\usepackage{booktabs}       % professional-quality tables
\usepackage{amsfonts}       % blackboard math symbols
\usepackage{nicefrac}       % compact symbols for 1/2, etc.
\usepackage{microtype}      % microtypography
\usepackage{lipsum}
\usepackage{fancyhdr}       % header
\usepackage{graphicx}       % graphics
\graphicspath{{media/}}     % organize your images and other figures under media/ folder
\usepackage{float}
\usepackage{framed,multirow}
\usepackage{amsmath}
\allowdisplaybreaks[4]
\usepackage{amssymb}
\usepackage{latexsym}
\usepackage{subcaption}
\usepackage{geometry}
\usepackage{pifont}  
\usepackage{xcolor}      
\usepackage{colortbl}     
\usepackage{bm}
\newtheorem{remark}{Remark}
\usepackage{algorithm}
\usepackage{algorithmic}
\usepackage{adjustbox}
\newcommand{\cmark}{\ding{51}}
\newcommand{\xmark}{\ding{55}}
\usepackage{tabularx}
\usepackage{pdflscape}
\newcolumntype{C}{>{\centering\arraybackslash}X}

\title{A phase-field neural solver for moving contact line problems with dynamic boundary conditions}

\author{
  Ziyan Chen \\
  Department of Mathematics, University of Macau \\
  Macao SAR, China \\
  \And
  Jinpeng Zhang \\
  Department of Mathematics, University of Macau \\
  Macao SAR, China \\
  \And
  Pai Zhang \\
  Department of Mathematics, University of Macau \\
  Macao SAR, China \\
  \And
  Li Luo \thanks{Corresponding authors} \\
  Department of Mathematics, University of Macau \\
  Macao SAR, China \\
  \texttt{liluo@um.edu.mo} \\
}

\begin{document}
\maketitle

\begin{abstract}
Phase-field models based on the Cahn--Hilliard equation coupled with dynamic boundary conditions provide a thermodynamically consistent framework for moving contact line (MCL) problems. Although physics-informed neural networks (PINNs) offer a mesh-free approach for solving partial differential equations, their direct application to MCL problems remains challenging due to long-time error accumulation, sharp interfacial profiles, localized contact line dynamics, and complex contact angle evolution. In this work, we propose MCL-PINNs, a specialized phase-field neural solver designed for MCL problems with dynamic boundary conditions. The method is built on a discrete-time formulation and incorporates several key techniques, including a multi-network time-marching scheme, a relaxed distribution constraint on the neural network outputs, variable scaling for sharply varying solution features, adaptive loss weighting, adaptive collocation sampling with interface extraction, and, when applicable, symmetry preservation through neural network inputs. These techniques improve the capability of the neural solver in resolving sharp interfacial profiles and contact line motion. The proposed method is validated through three numerical examples involving droplet coalescence, shear-induced droplet deformation, and dynamic wetting in a heterogeneous channel. The numerical results show that MCL-PINNs significantly improve prediction accuracy and robustness compared with standard PINNs formulations, enabling reliable resolution of complex interfacial evolution and moving contact line dynamics.
\end{abstract}

% keywords can be removed

\section{Introduction} \label{sec:intro}
Phase field modeling, originating in the classical formulation by Cahn and Hilliard \cite{cahn1958free}, has established itself as a robust framework for simulating interfacial evolution and multiphase dynamics. This approach has been applied to a broad range of complex physical systems, including two-phase flows \cite{Anderson1998Diffuse,zhang2026relaxed,zhang2025fully}, material phase transitions \cite{Chen2002Phase}, fracture mechanics \cite{Miehe2010Thermo}, and biomembrane evolution \cite{du2006simulating}. 
By describing the interface as a diffuse region with a continuous order parameter, the phase-field method can naturally handle changes in interface topology. Therefore, it avoids the need for explicit interface reconstruction and remeshing, which are often required in conventional sharp-interface tracking methods \cite{Lowengrub1998Quasi}.
For problems involving solid-liquid-gas three-phase contact, dynamic boundary conditions are needed to describe moving contact line (MCL) phenomena, such as wetting and spreading, the evolution of the dynamic contact angle, and contact angle hysteresis (CAH) \cite{qian2006variational, xu2010wenzel_cassie, xu2011hysteresis, wang2016dynamic}. 

Phase-field models for MCL problems, typically based on the fourth-order Cahn--Hilliard equation, are still difficult to solve numerically. When these models are discretized by traditional methods, such as finite element or finite difference methods, the numerical results may be highly sensitive to the mesh resolution, especially when the motion of a thin diffuse interface must be accurately resolved. In such cases, sufficiently fine meshes are required near the interface and the contact line, which can substantially increase the computational cost in large-scale or high-dimensional problems. Moreover, the implementation of traditional solvers is often technically demanding, as it requires the careful discretization of high-order nonlinear terms, the design of stable time-integration schemes, the treatment of coupled bulk-boundary dynamics, and the accurate enforcement of complicated dynamic boundary conditions \cite{luo2017parallel,luo2017efficient,luo2025numerical}. These challenges may lead to strong mesh dependence and increased algorithmic complexity, while limiting the adaptability of traditional solvers to evolving interfaces and complex geometries, especially in problems involving moving contact lines and heterogeneous boundaries. They therefore motivate the development of alternative computational paradigms that reduce the reliance on mesh-based discretization.

Recent advances in computational hardware and optimization algorithms have laid the computational foundation for the practical implementation of Physics-Informed Neural Networks (PINNs) \cite{raissi2019physics}. By embedding governing partial differential equations (PDEs), along with initial and boundary conditions, directly into the neural network loss functional, PINNs provide a unified, mesh-free framework for solving both forward and inverse problems. This method has demonstrated notable efficacy in parameter identification and complex system modeling \cite{raissi2020hidden, wang2021understanding}. Systematic reviews \cite{cuomo2022scientific, tanyu2023deep}, along with benchmark suites such as PINNACLE \cite{zhongkai2024pinnacle}, have established PINNs as a powerful computational tool for scientific computing.

PINNs have been used to solve phase-field models for interfacial evolution problems \cite{wight2020solving, qiu2022physics, chen2024pf, chen2025sharp, qiu2025direct}. In these studies, the Cahn--Hilliard equation is usually included in the loss function of the neural network, showing the potential of PINN-type methods for diffuse-interface problems. However, most existing works mainly focus on relatively standard phase-field settings, such as bulk phase separation, interface motion, or simplified boundary conditions. The moving contact line problem considered here is more difficult because it involves the coupled evolution of sharp interfacial layers, moving contact lines, and relaxing contact angles.
These difficulties are not fully addressed by existing PINN-based phase-field solvers. First, the sharp phase-field interface is difficult to represent with standard neural network outputs. Since the order parameter changes rapidly between $-1$ and $1$ across a thin diffuse layer, the network may struggle to capture this transition accurately and stably. Second, phase-field evolution over long time intervals can suffer from error accumulation, so a single network trained over the whole time domain may not give reliable long-time predictions \cite{mattey2022novel}. Third, the important features of MCL dynamics are concentrated near the diffuse interface and the contact line. Standard collocation sampling may miss these localized regions or require a very large number of points, which can reduce both accuracy and efficiency \cite{wu2023comprehensive, krishnapriyan2021characterizing}. Finally, dynamic wetting involves nonlinear and spatially varying boundary conditions. If these boundary conditions are imposed only through fixed soft penalty terms, the predicted boundary behavior may be inaccurate, especially for contact angle hysteresis and contact line pinning \cite{yue2020thermodynamically, liu2022unified}.

Although existing PINN-based phase-field studies provide useful ideas for diffuse-interface modeling, they are not directly sufficient for the MCL problems studied in this work. To address this difficulty, we develop MCL-PINNs, a unified computational framework specifically designed for the Cahn--Hilliard equation with dynamic contact angle boundary conditions. By incorporating several key techniques tailored to moving contact line dynamics, the proposed method effectively captures the coupled interfacial and boundary evolution while improving numerical stability and solution accuracy. Importantly, instead of using standard continuous-time formulations that rely on automatic differentiation of temporal derivatives, MCL-PINNs is built upon the \textit{discrete-time physics-informed neural networks} paradigm \cite{raissi2019physics}, with a neural architecture specifically designed to incorporate the implicit Runge--Kutta discretization. In this formulation, the temporal evolution is represented through coupled stage-value equations, thereby reformulating the evolutionary governing equation as a sequence of spatial boundary-value problems evaluated at discrete implicit Runge--Kutta nodes.
 This discrete-time formulation helps reduce the influence of stiff gradient accumulation and temporal error propagation often observed in continuous-time implementations. 
 
To facilitate long-time simulations, the framework employs a \textit{time-marching strategy} with multiple networks, in which the full time interval is divided into a sequence of shorter temporal windows trained in a marching manner. This design helps reduce error accumulation across long temporal horizons \cite{jagtap2020extended,meng2020ppinn,qiu2022physics,mattey2022novel}. Within each discrete temporal window, a \textit{variable-scaling} transformation is applied to rescale the governing system, spreading sharply varying solution features over a broader coordinate range and reducing their effective frequency for neural-network approximation \cite{ko2025vs}. In addition, physical and geometric properties are also incorporated through tailored architectural designs. We introduce a \textit{relaxed distribution constraint}, which employs a scaled latent mapping to bound the required latent activation range, thereby avoiding excessively large bulk activations while maintaining the approximation capability for sharp interfacial layers \cite{rahaman2019spectral, lu2021physics}. When the problem configuration is inherently symmetric, a \textit{symmetry-preserving} input transformation is further employed to embed the symmetry into the neural network representation and restrict the approximation to a symmetry-consistent function space. To better resolve localized interfacial structures without relying on excessively dense global sampling, an \textit{adaptive sampling strategy} is devised to identify regions near the diffuse interface and moving contact line, and to adjust the collocation-point distribution accordingly \cite{hou2023enhancing,wu2023comprehensive}. Moreover, to balance the different components of the composite loss, an \textit{adaptive loss weighting} scheme is further employed to adjust the relative contributions of the PDEs residuals and MCL boundary terms, alleviating training imbalance caused by dominant loss components \cite{wang2024respecting}. Numerical experiments demonstrate that MCL-PINNs significantly improves prediction accuracy and training robustness compared with standard PINNs formulations \cite{raissi2019physics}, enabling accurate simulations of complex wetting and multiphase interfacial dynamics, including contact line motion and dynamic contact angle hysteresis.

The remainder of this paper is organized as follows. Section~\ref{sec:model} introduces the mathematical formulation of the phase-field model for moving contact line dynamics. Section~\ref{sec:methodology} presents the MCL-PINNs framework, with emphasis on the key techniques and their integration into the overall neural solver. Section~\ref{sec:experimental_results} validates the proposed method through three numerical examples, including droplet coalescence, shear-induced droplet deformation, and dynamic wetting through a heterogeneous channel. Finally, Section~\ref{sec:conclusion} provides concluding remarks.

\section{Modeling of the moving contact line problems} \label{sec:model}
Let $\Omega \subset \mathbb{R}^d$ ($d=2,3$) be a bounded domain. 
We consider the interfacial dynamics of a binary mixture governed by the Cahn--Hilliard equation which can be written as follows,
\begin{equation} \label{eq:CH_PDE1}
    \left\{
    \begin{aligned}
        & \frac{\partial \phi}{\partial t} = M\Delta\mu, && \text{in } \Omega,\\
        & \mu= \frac{\delta \mathcal{F}_{bulk}}{\delta \phi} = -\varepsilon^2\Delta\phi + F'(\phi), && \text{in } \Omega.
    \end{aligned}
    \right.
\end{equation}
This model provides a continuum description for the evolution of a conserved order parameter, $\phi(\mathbf{x}, t)$, which distinguishes between two immiscible phases. $\mu$ is the chemical potential, and $M$ represents the mobility. $\varepsilon$ is the parameter which is proportional to the interface thickness. The formulation is rooted in a Ginzburg-Landau free energy functional, $\mathcal{F}_{bulk}(\phi)$, that accounts for both the bulk free energy of the mixture and the energetic penalty associated with spatial gradients of the order parameter. The bulk free energy is given by
\begin{equation}
    \mathcal{F}_{bulk}(\phi) = \int_{\Omega} \left( \frac{\varepsilon^2}{2} |\nabla \phi|^2 + F(\phi) \right)\, d\mathbf{x},
\end{equation}
where $F(\phi) = \frac{1}{4}(\phi^2 - 1)^2$ is the double-well potential ensuring phase separation into pure components corresponding to $\phi = \pm 1$. 
\begin{figure}
    \centering
    \includegraphics[width=0.8\textwidth]{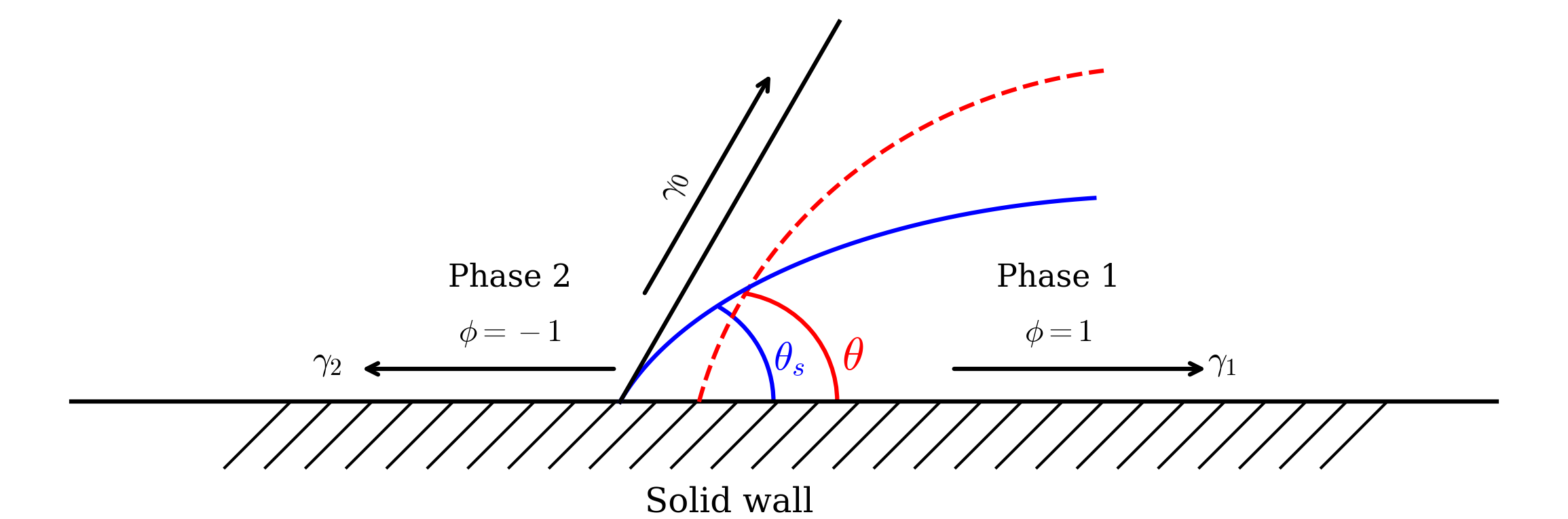}
    \caption{Schematic illustration of interfacial tension balance and contact angle definition at a solid wall in a two-phase system. The blue solid and red dashed curves denote the interface profiles corresponding to the static contact angle $\theta_s$ and the dynamic contact angle $\theta$, respectively. $\gamma_{0}$, $\gamma_{1}$, and $\gamma_{2}$ denote the interfacial tensions of Phase~1-Phase~2, Phase~1-Solid wall, and Phase~2-Solid wall, respectively.}
    \label{fig:contact_angle}
\end{figure}

When the interface intersects the boundary, the problem becomes a moving contact line problem. At equilibrium, the angle between the fluid-fluid interface and the solid wall is known as the static contact angle, denoted by $\theta_s$. 
In non-equilibrium situations when  the contact line moves along the solid boundary, the contact angle becomes dynamic and depends not only on the local interfacial geometry but also on additional energy dissipation due to microscopic slip or molecular reorganization near the solid surface. To account for this effect, we consider the interfacial free energy per unit area at the fluid-solid interface
\begin{equation}
\gamma(\phi)=-\frac{\sqrt{2}}{3}\cos\theta_s \sin\!\left(\frac{\pi}{2}\phi\right).
\end{equation}
The corresponding surface energy at the fluid-solid interface is $\mathcal{F}_{surf}(\phi)=\varepsilon\int_{\partial \Omega} \gamma(\phi)\, ds$.
Applying the Onsager variational principle \cite{qian2006variational} to 
the Cahn--Hilliard system with boundary contributions yields the following dynamic boundary conditions,
\begin{align} 
    & \varepsilon \frac{\partial \phi}{\partial \mathbf{n}} = -\frac{\partial\gamma(\phi)}{\partial\phi} - \frac{1}{\alpha}\frac{\partial \phi}{\partial t}, \quad \text{ on } \partial \Omega, \label{eq:relaxation_bc}\\
    & \nabla \mu \cdot \mathbf{n} = 0, \quad \text{ on } \partial \Omega, \label{eq:conservation_bc}
\end{align}
where $\mathbf{n}$ denotes the outward unit normal vector on $\partial \Omega$. 
The first condition \eqref{eq:relaxation_bc} is a relaxation-type boundary condition, where $\alpha$ denotes the relaxation parameter. It accounts for energy dissipation at the boundary and is capable of capturing the motion of moving contact lines as well as the dynamic behavior of wetting phenomena.
In the equilibrium state, the Young’s equation for the static contact angle $\theta_s$ can be derived from the boundary condition~\eqref{eq:relaxation_bc},
\begin{equation}
\gamma_{0}\cos\theta_s = \gamma_2 - \gamma_1,
\end{equation}
where $\gamma_{0}$, $\gamma_{1}$, and $\gamma_{2}$ denote the interfacial tensions of Phase~1-Phase~2, Phase~1-Solid wall, and Phase~2-Solid wall, respectively. See \cite{qian2006variational} more details. 
Figure~\ref{fig:contact_angle} shows a schematic illustration of interfacial tension balance and contact angle definition at a solid wall. The second condition \eqref{eq:conservation_bc} ensures that the total mass within the domain is conserved by preventing any mass flux across the boundary. 

The Cahn--Hilliard equation \eqref{eq:CH_PDE1}, equipped with the dynamic boundary conditions \eqref{eq:relaxation_bc}-\eqref{eq:conservation_bc}, remains a gradient flow associated with the following total free energy
\begin{equation}
\mathcal{F}(\phi)=\mathcal{F}_{bulk}(\phi)+\mathcal{F}_{surf}(\phi),
\end{equation}
and satisfies the following energy
dissipation law
\begin{equation}
\frac{d}{dt} \mathcal{F}(\phi)
= - M\int_{\Omega} \left| \nabla \mu \right|^2 \,d\mathbf{x}-\varepsilon\alpha\int_{\partial \Omega}\left|L(\phi)\right|^2\,ds,
\end{equation}
where 
\begin{equation}
L(\phi)=\left.\frac{\delta\mathcal{F}(\phi)}{\delta\phi}\right|_{\partial\Omega}=\varepsilon \frac{\partial \phi}{\partial \mathbf{n}}+\frac{\partial\gamma\left(\phi\right)}{\partial\phi},
\end{equation}
yielding a thermodynamically consistent mechanism for contact line motion.

The initial condition for $\phi$ is given by
\begin{equation} \label{eq:CH_initial}
\phi(x, y, 0) = \phi^0(x, y).
\end{equation}

\begin{remark} \label{remark1}
    When the interface moves under the action of an external force, an advection term, $\mathbf{u}\cdot\nabla\phi$, is added to the left-hand side of \eqref{eq:CH_PDE1} to describe this effect, where $\mathbf{u}$ is the velocity field driving the interface motion.
\end{remark}

\section{The phase-field neural solver} \label{sec:methodology}
This section presents the proposed framework in detail for solving the Cahn--Hilliard equation subject to dynamic boundary conditions.

\subsection{Time-marching scheme with multiple networks}
\label{subsec:time_marching_scheme}
For the Cahn--Hilliard equation with dynamic boundary conditions, training a single physics-informed neural network over the entire spacetime domain often leads to slow convergence or error accumulation, primarily due to temporal causality violations and the stiff, multiscale nature of the MCL dynamics.
To address these issues, we adopt a \emph{time-marching strategy} \cite{wight2020solving} that decomposes the global problem into a sequence of local-in-time subproblems.

\begin{figure}
    \centering
    \includegraphics[width=0.95\linewidth]{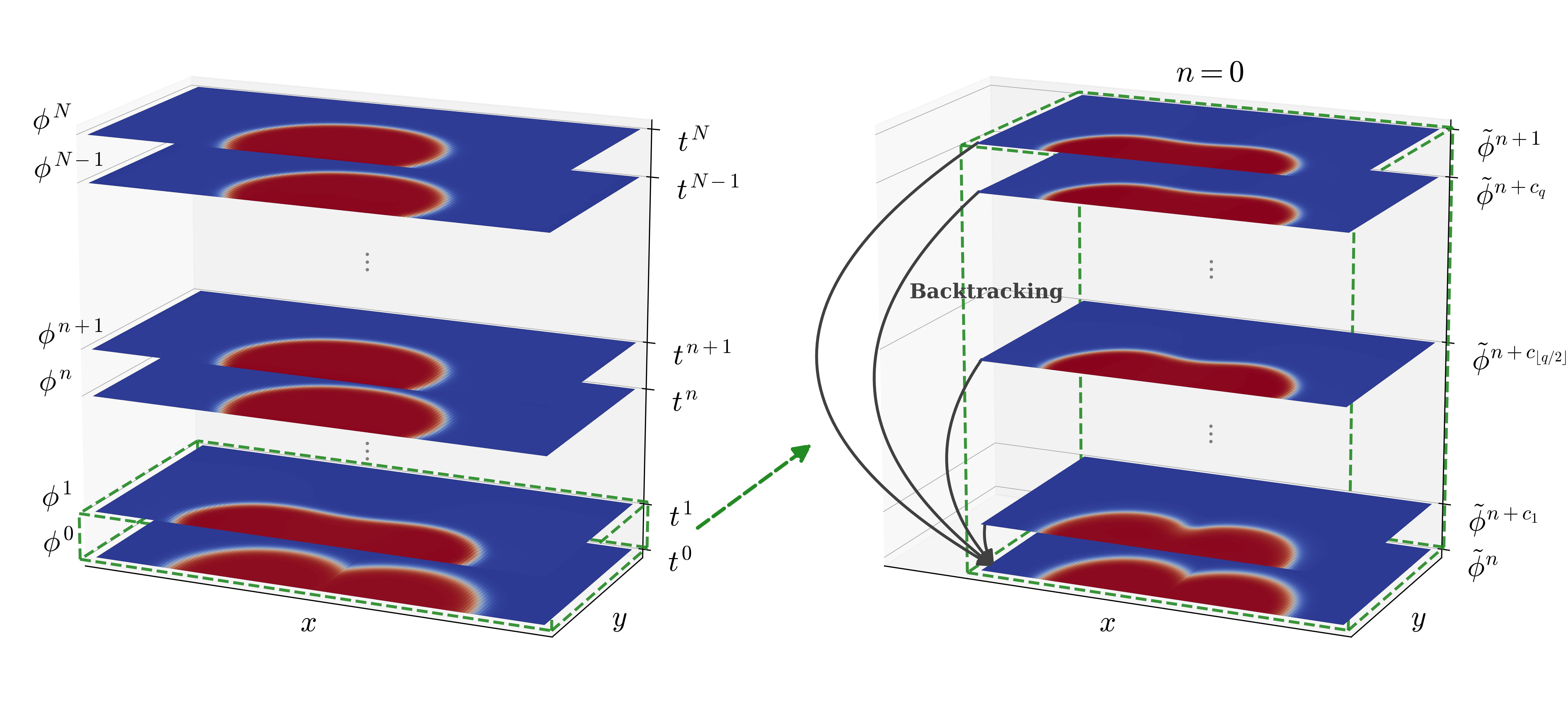}
    \caption{A time-marching framework embedded with multiple Implicit Runge--Kutta (IRK) discretized PINNs. 
(Left) Global view of the sequential time-marching scheme over discrete time levels $ \{t^0, t^1, \dots, t^n, t^{n+1} \dots, t^{N-1}, t^N\} $; 
(Right) Internal architecture of the PINNs associated with a single time interval  $[t^{n}, t^{n+1}]$, where the IRK stages are encoded as intermediate outputs of a shared neural network and coupled through a backtracking formulation to construct the residual of the Cahn--Hilliard equation.
}
    \label{fig:time_marching_figure}
\end{figure}

Specifically, for a given final time $T > 0$, the interval $[0, T]$ is partitioned into $N$ contiguous, non-overlapping intervals 
\begin{equation}
    [t^0, t^1],\, [t^1, t^2],\, \dots,\, [t^{N-1}, t^N],
\end{equation}
where $t^0 = 0$ and $t^N = T$.
A separate neural network is employed to approximate the solution on each interval $[t^{n}, t^{n+1}]$, $n = 0, \dots, N-1$.
The training proceeds sequentially. 
The first network is trained on $\Omega \times [t^0, t^1]$ using the given initial condition at $t = t^0$. 
Upon convergence, its predicted solution at $t = t^1$ is used as the initial condition for the second network on $[t^1, t^2]$.  To accelerate convergence, each subsequent network inherits the trained parameters of its predecessor as a warm start, since the solution at adjacent time windows is expected to be similar.
This procedure is carried out iteratively, with the approximation generated by network~$n$ at time~$t^n$ acting as the initial condition for network~$n+1$.
The solution evolves through time across consecutive intervals. The sequential training procedure is illustrated in Figure~\ref{fig:time_marching_figure}.
This approach reduces the effective learning complexity for each individual network, thereby promoting more stable and efficient training. 

\subsection{Discrete-time physics-informed neural networks} \label{subsec:discrete_time_model_PINNs}
To address the temporal evolution of the Cahn--Hilliard system in each time interval $[t^{n}, t^{n+1}]$, we employ a \emph{discrete-time} PINNs framework. This approach discretizes the governing equations in time using a $q$-stage Implicit Runge--Kutta (IRK) method \cite{iserles2009first}, reducing the evolutionary equation to a system of coupled boundary-value problems at each time step. 

Consider $\Omega \subset \mathbb{R}^2$ and $\mathbf{x}=(x,y)$.
Let the Cahn--Hilliard system be written in the general form $\partial_t \phi = \mathcal{M}\left[\phi, \mu\right]$, where the nonlinear operator is given by $\mathcal{M}\left[\phi, \mu\right] = -\mathbf{u}\cdot\nabla \phi + M\Delta \mu$ with $\mu = -\varepsilon^2 \Delta \phi + F'(\phi) = -\varepsilon^2 \Delta \phi - \phi + \phi^3$.  Here $\mathbf{u}$ denotes an external velocity field, which is present only when advection is considered (cf.\ the remark in Section~\ref{sec:model}). Applying a $q$-stage Runge--Kutta scheme with a uniform time step $\Delta t = T/N$ to advance the solution from time $t^n$ to $t^{n+1} = t^n + \Delta t$, we obtain the stage values $\phi^{n+c_i}$ and the final solution $\phi^{n+1}$:
\begin{equation} \label{eq:ch_runge_kutta}
\begin{aligned}
    \phi^{n+c_i} &= \phi^n + \Delta t \sum_{j=1}^{q} a_{ij} \mathcal{M}\left[\phi^{n+c_j},\mu^{n+c_j}\right], \quad i=1, \ldots, q, \\
    \phi^{n+1} &= \phi^n + \Delta t \sum_{j=1}^{q} b_{j} \mathcal{M}\left[\phi^{n+c_j},\mu^{n+c_j}\right].
\end{aligned}
\end{equation}
Here, $\phi^{n+c_j}(x, y) = \phi(x, y, t^n + c_j \Delta t)$ for $j = 1, \dots, q$. We adopt the convention $c_0 = 0$ and $c_{q+1} = 1$, such that $\phi^{n+c_0} = \phi^n$ denotes the known state at the beginning of the time interval, while $\phi^{n+c_{q+1}} = \phi^{n+1}$ represents the target state to be predicted by the discrete-time PINNs. This unified formulation encompasses both explicit and implicit time-integration schemes, determined by the specific selection of the Runge--Kutta coefficients $\{a_{ij}, b_j, c_j\}$. In this work, we focus on IRK schemes to ensure unconditional stability for stiff phase field dynamics.

Equations \eqref{eq:ch_runge_kutta} can be equivalently expressed as
\begin{equation}
\begin{array}{l}
\phi^{n}=\phi_{i}^{n}, \quad i=1, \ldots, q, \\
\phi^{n}=\phi_{q+1}^{n},
\end{array}
\end{equation}
where
\begin{equation} \label{eq:ch_operator_form}
\begin{aligned}
    \phi_i^n &:= \phi^{n+c_i} - \Delta t \sum_{j=1}^{q} a_{ij} \mathcal{M}\left[\phi^{n+c_j},\mu^{n+c_j}\right], \quad i=1, \ldots, q, \\
    \phi_{q+1}^n &:= \phi^{n+1} - \Delta t \sum_{j=1}^{q} b_{j} \mathcal{M}\left[\phi^{n+c_j},\mu^{n+c_j}\right].
\end{aligned}
\end{equation}

We begin by assuming that the quantities
\begin{equation} \label{eq:nn_output}
    \left[ \phi^{n+c_1}(x, y), \ldots, \phi^{n+c_q}(x, y), \phi^{n+1}(x, y),  \mu^{n+c_1}(x, y), \ldots, \mu^{n+c_q}(x, y), \mu^{n+1}(x, y) \right],
\end{equation}
are jointly approximated by the outputs of a neural network
\begin{equation} \label{eq:nn_output_approx}
    \left[\tilde{\phi}^{n+c_1}(x, y), \ldots, \tilde{\phi}^{n+c_q}(x, y), \tilde{\phi}^{n+1}(x, y),  \tilde{\mu}^{n+c_1}(x, y), \ldots, \tilde{\mu}^{n+c_q}(x, y), \tilde{\mu}^{n+1}(x, y) \right].
\end{equation}
Substituting the neural network approximations \eqref{eq:nn_output_approx} into the back-projection formula \eqref{eq:ch_operator_form}, we obtain the quantities $\tilde{\phi}_j^n$ for $j=1,\ldots,q+1$. This yields the prediction
\begin{equation} \label{eq:nn_output2}
    \left[\tilde{\phi}^{n}_1(x, y), \ldots, \tilde{\phi}^{n}_q (x, y), \tilde{\phi}^{n}_{q+1}(x, y),  \tilde{\mu}^{n+c_1}(x, y), \ldots, \tilde{\mu}^{n+c_q}(x, y), \tilde{\mu}^{n+1}(x, y) \right].
\end{equation}

The network parameters $\boldsymbol{\theta}$ are then determined by minimizing the discrepancy between these back-projected values and the known state $\phi^n$, together with the residuals of the boundary conditions and the chemical potential definition, which are grouped into the following composite loss function
\begin{equation} \label{eq:ch_total_loss}
    \mathcal{L}(\boldsymbol{\theta}) =  
    \lambda_{r_{\phi}} \mathcal{L}_{r_{\phi}}(\boldsymbol{\theta}) + \lambda_{r_{\mu}} \mathcal{L}_{r_{\mu}}(\boldsymbol{\theta}) + \lambda_{b_{\phi}} \mathcal{L}_{b_{\phi}}(\boldsymbol{\theta}) + \lambda_{b_{\mu}} \mathcal{L}_{b_{\mu}}(\boldsymbol{\theta}),
\end{equation}
where $\lambda_{r_\phi}$, $\lambda_{r_\mu}$, $\lambda_{b_\phi}$, and $\lambda_{b_\mu}$ are self-adaptive weighting parameters to be introduced in Section~\ref{subsec:Adaptive_loss_weighting}. The residual losses $\mathcal{L}_{r_{\phi}}$ and $\mathcal{L}_{r_{\mu}}$ enforce the condition that the stage-wise predictions of the neural network satisfy the IRK discretization of the Cahn--Hilliard equation. Specifically, $\mathcal{L}_{r_{\phi}}$ represents the discrepancy between the network predictions projected back to time level $t^n$ (as defined in \eqref{eq:ch_operator_form}) and the known solution $\phi^n$ at $N_r$ interior collocation points $\{(x_r^k, y_r^k)\}_{k=1}^{N_r}$, 
\begin{equation} \label{eq:ch_loss_f_phi}
\mathcal{L}_{r_{\phi}}(\boldsymbol{\theta}) = \frac{1}{(q+1)N_r} \sum_{j=1}^{q+1} \sum_{k=1}^{N_r}  \left| \tilde{\phi}_j^n(x_r^k, y_r^k; \boldsymbol{\theta}) - \phi^n(x_r^k, y_r^k)\right|^2.
\end{equation}
The term $\mathcal{L}_{r_{\mu}}$ enforces the definition of the chemical potential,
\begin{equation} \label{eq:ch_loss_f_mu}
\mathcal{L}_{r_{\mu}}(\boldsymbol{\theta}) = \frac{1}{qN_r} \sum_{j=1}^{q} \sum_{k=1}^{N_r}
\left|\mathcal{R}\left[\tilde{\phi}^{n+c_j}(x_r^k, y_r^k; \boldsymbol{\theta}), \tilde{\mu}^{n+c_j}(x_r^k, y_r^k; \boldsymbol{\theta})\right] \right|^2,
\end{equation}
where
\begin{align}\label{eq:loss_R}
\mathcal{R}\left[\tilde{\phi}^{n+c_j}(x_r^k, y_r^k; \boldsymbol{\theta}), \tilde{\mu}^{n+c_j}(x_r^k, y_r^k; \boldsymbol{\theta})\right] &= \tilde{\mu}^{n+c_j}(x_r^k, y_r^k; \boldsymbol{\theta}) \nonumber \\
&- \left(-\varepsilon^2 \Delta \tilde{\phi}^{n+c_j}(x_r^k, y_r^k; \boldsymbol{\theta}) - \tilde{\phi}^{n+c_j}(x_r^k, y_r^k; \boldsymbol{\theta}) + ({\tilde{\phi}^{n+c_j}})^3(x_r^k, y_r^k; \boldsymbol{\theta})\right).
\end{align}

For the initial time interval ($n=0$), $\phi^0(x, y)$ is given by the initial condition. For subsequent time intervals ($n>0$), $\phi^n(x, y)$ is the solution predicted by the network trained at the previous time step, as described in Section~\ref{subsec:time_marching_scheme}.

\begin{remark}
Note that the back-projection residual in \eqref{eq:ch_loss_f_phi} is enforced at all $q+1$ time levels (including $t^{n+1}$), whereas the chemical potential constitutive residual in \eqref{eq:ch_loss_f_mu} is enforced only at the $q$ IRK stage points $t^{n+c_1}, \ldots, t^{n+c_q}$. This is because the IRK update formula~\eqref{eq:ch_operator_form} involves $\mu^{n+c_j}$ only for $j = 1, \ldots, q$; the chemical potential at $t^{n+1}$ does not enter the stage equations and is therefore constrained indirectly through the boundary loss.
\end{remark}

The boundary losses $\mathcal{L}_{b_{\phi}}(\boldsymbol{\theta})$ and $\mathcal{L}_{b_{\mu}}(\boldsymbol{\theta})$ penalize deviations from the specified boundary conditions at a set of $N_b$ boundary points $\{(x_b^k, y_b^k)\}_{k=1}^{N_b}$ for all $q+1$ IRK stages, specifically,
\begin{equation} \label{eq:ch_loss_b_phi}
\mathcal{L}_{b_{\phi}}(\boldsymbol{\theta}) = \frac{1}{(q+1)N_b} \sum_{j=1}^{q+1}\sum_{k=1}^{N_b}  \left| \mathcal{B_{\phi}}\left[\tilde{\phi}^{n+c_j}(x_b^k, y_b^k;\boldsymbol{\theta})\right]\right|^2 ,
\end{equation}
and
\begin{equation} \label{eq:ch_loss_b_mu}
\mathcal{L}_{b_{\mu}}(\boldsymbol{\theta}) = \frac{1}{(q+1)N_b} \sum_{j=1}^{q+1}\sum_{k=1}^{N_b} \left| \mathcal{B_{\mu}}\left[\tilde{\mu}^{n+c_j}(x_b^k, y_b^k;\boldsymbol{\theta})\right] \right|^2 .
\end{equation}
Here, the boundary operators $\mathcal{B_{\phi}}$ and $\mathcal{B_{\mu}}$ are given by
\begin{equation} \label{eq:loss_Bphi}
\mathcal{B_{\phi}}\left[\tilde{\phi}^{n+c_j}(x_b^k, y_b^k;\boldsymbol{\theta})\right]=
\varepsilon \frac{\partial \tilde{\phi}^{n+c_j}(x_b^k, y_b^k;\boldsymbol{\theta})}{\partial \mathbf{n}} - \frac{\sqrt{2} \pi}{6} \cos \theta_s \cos\left(\frac{\pi}{2}\tilde{\phi}^{n+c_j}(x_b^k, y_b^k;\boldsymbol{\theta})\right) + \frac{1}{\alpha}\frac{\partial \tilde{\phi}^{n+c_j}(x_b^k, y_b^k;\boldsymbol{\theta})}{\partial t},
\end{equation}
and
\begin{equation}
\label{eq:loss_Bmu}
\mathcal{B_{\mu}}\left[\tilde{\mu}^{n+c_j}(x_b^k, y_b^k;\boldsymbol{\theta})\right]
    = \nabla \tilde{\mu}^{n+c_j}(x_b^k, y_b^k;\boldsymbol{\theta}) \cdot \mathbf{n}.
\end{equation}

Define $t^{n+c_j}=t^n+c_j\Delta t$ for $j=0,1,\ldots,q+1$. 
Then, for $j=1,2,\ldots,q+1$, the time derivative 
$\partial_t \tilde{\phi}^{n+c_j}$ at the time level $t^{n+c_j}$ is approximated by the backward finite difference
\begin{equation} \label{eq:FDM_phi_t}
    \frac{\partial \tilde{\phi}^{n + c_j}(x_b^k, y_b^k;\boldsymbol{\theta})}{\partial t}\approx 
    \frac{\tilde{\phi}^{n + c_j}(x_b^k, y_b^k;\boldsymbol{\theta}) - \tilde{\phi}^{n + c_{j-1}}(x_b^k, y_b^k;\boldsymbol{\theta})}{t^{n + c_j} - t^{n + c_{j-1}}}, \quad j = 1, 2, \dots, q+1.
\end{equation}
The spatial derivatives $\partial_{\mathbf{n}}\tilde{\phi}^{n+c_j}$ in \eqref{eq:loss_Bphi} are obtained via automatic differentiation from the neural network output~\eqref{eq:nn_output_approx}.

\begin{figure}
    \centering
\includegraphics[width=\textwidth]{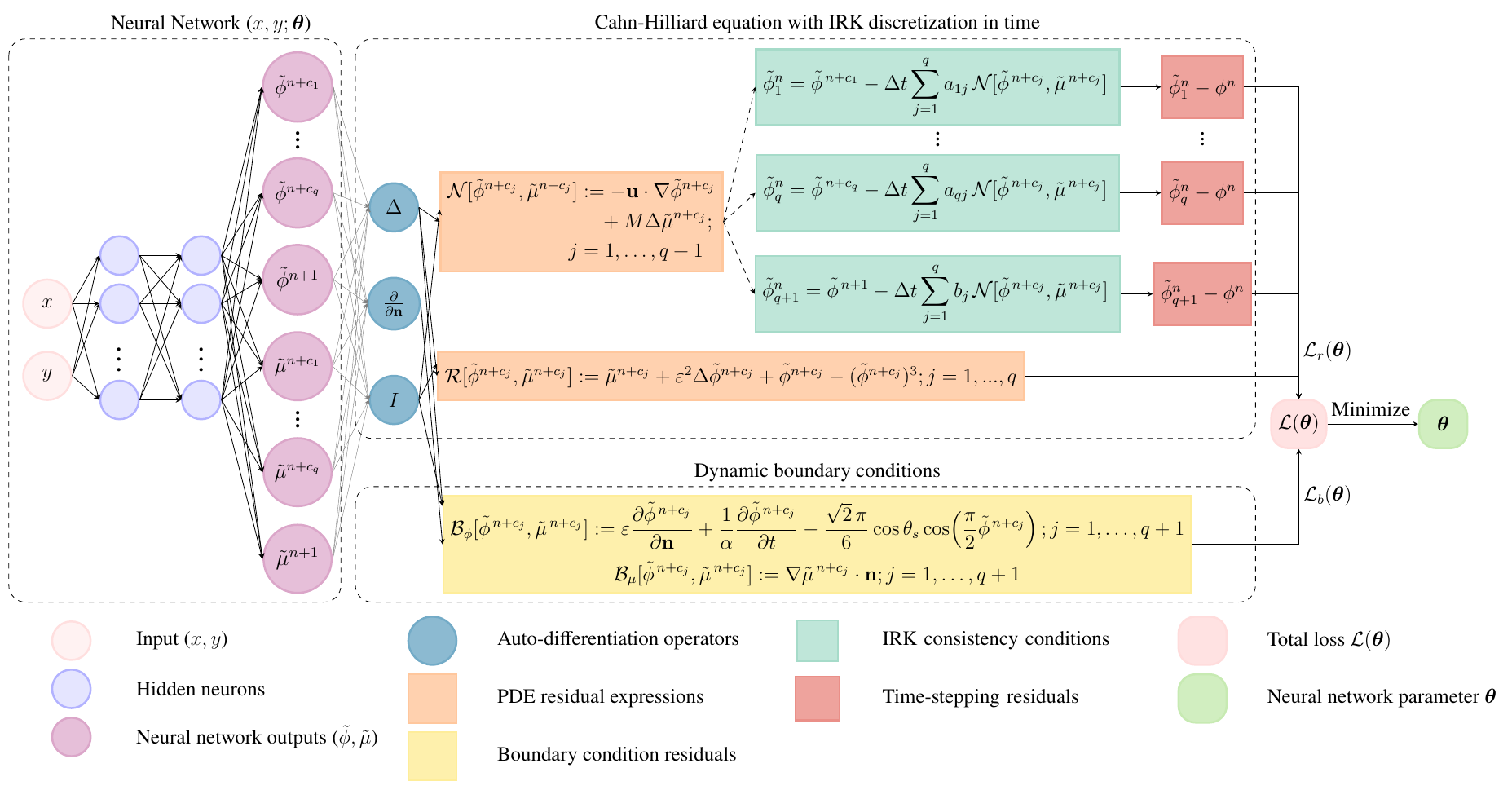}
    \caption{Discrete-time physics-informed neural network framework based on Implicit Runge--Kutta (IRK) method for solving the moving contact line problem modeled by the Cahn--Hilliard equation with dynamic wetting boundary conditions.}
    \label{fig:discrete_time_model}
\end{figure}

Figure~\ref{fig:discrete_time_model} presents a schematic illustration of this discrete-time PINNs framework. The architecture comprises a fully connected neural network mapping $(x,y)$ to the IRK stage values $\tilde{\phi}^{n+c_j}, \tilde{\mu}^{n+c_j}$, an automatic differentiation module for computing spatial derivatives, and residual blocks encoding the governing equations and boundary conditions. Stage-wise consistency is enforced by regressing the back-projected values $\tilde{\phi}_j^n$ against the known state $\phi^n$. 

\subsection{Relaxed distribution constraint}
\label{subsec:RDC}

\begin{figure}
    \centering
    \includegraphics[width=\linewidth]{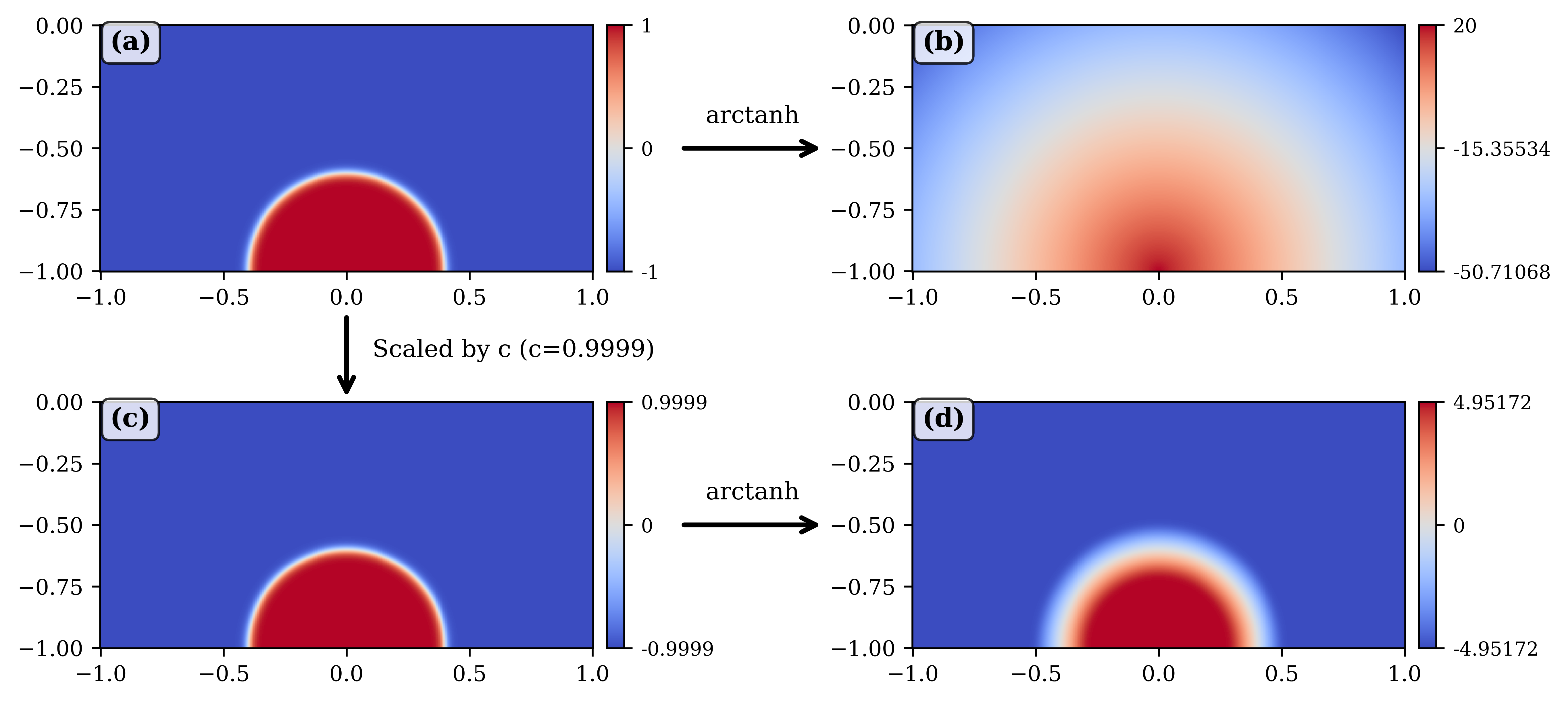}
    \caption{Illustration of the relaxed distribution constraint (RDC) mapping. 
(a) Target phase-field variable $\phi$, which exhibits a sharp diffuse-interface transition between the bulk phases. 
(b) Network output variable $z=\operatorname{arctanh}(\phi)$ required by the standard distribution constraint, where the interfacial jump is mapped into a smooth transition but the bulk phases require a large output range. 
(c) Scaled phase-field variable $c\phi$ with $c=0.9999$. 
(d) Relaxed output variable $z=\operatorname{arctanh}(c\phi)$, for which the scaling factor bounds the required output range while retaining the smooth transition structure. 
    }
    \label{fig:rdc}
\end{figure}

Solutions of the Cahn--Hilliard equation typically exhibit a sharp interfacial profile: the order parameter $\phi$ rapidly transitions between the two bulk values $\pm 1$ across a diffuse interface of width $\mathcal{O}(\varepsilon)$. Such localized transitions are difficult to approximate using standard fully connected neural networks, whose smooth activation functions tend to favor low-frequency components and may therefore struggle to represent steep interfacial layers, a phenomenon commonly associated with spectral bias \cite{rahaman2019spectral}.

To incorporate this phase-field structure into the network representation, a natural strategy is to impose a \emph{distribution constraint} at the output level. Specifically, one may represent the phase variable as $\tilde{\phi}=\tanh(z)$, where $z$ denotes the output of the neural network. Unlike hard constraints designed to exactly enforce initial or boundary conditions \cite{lu2021physics}, this constraint targets the intrinsic shape of the solution. The $\tanh$ function provides a smooth monotone transition between $-1$ and $1$, which is consistent with the qualitative structure of the phase-field profile and reduces the burden on the hidden layers to learn the interfacial jump from scratch. However, this direct formulation also introduces a practical issue. Since $\tanh(z)$ approaches $\pm 1$ only in the limit $|z|\to\infty$, the network must generate large output values in the bulk regions where $\phi\approx\pm 1$, which can lead to ill-conditioning and slow convergence during training.

To address this issue, we introduce a \emph{relaxed distribution constraint} (RDC) by incorporating a scaling factor $c<1$:
\begin{equation}\label{eq:RDC}
    \tilde{\phi} = \frac{\tanh(z)}{c}.
\end{equation}
Consequently, the physical bulk values $\tilde{\phi}=\pm 1$ can be reached at finite output values, namely $z=\pm\operatorname{arctanh}(c)$. For instance, when $c=0.9999$, one has $\operatorname{arctanh}(c)\approx 4.95$. Therefore, the proposed constraint preserves the smooth transition structure induced by the $\tanh$ function while avoiding the need for excessively large outputs in the bulk phases. This relaxation improves the conditioning of the output representation and facilitates more stable optimization.

Figure~\ref{fig:rdc} illustrates the effect of this relaxation on the required output representation. Figure~\ref{fig:rdc}(a) shows a representative target phase field $\phi$ containing a circular droplet, spanning the full physical range $[-1,1]$. Under the standard distribution constraint ($c=1$), reproducing this field requires the network to output $z=\operatorname{arctanh}(\phi)$, displayed in Figure~\ref{fig:rdc}(b); since $\phi$ approaches $\pm 1$ in the bulk phases, the corresponding values of $z$ increase sharply and cover a wide range, approximately $[-50,20]$ in this example. 
Figures~\ref{fig:rdc}(c) and~\ref{fig:rdc}(d) show the relaxed mapping with $c=0.9999$. By first scaling the phase-field variable as $c\phi$, the required network output becomes $z=\operatorname{arctanh}(c\phi)$, which is confined to a much narrower range, approximately $[-4.95,4.95]$. Therefore, RDC substantially reduces the magnitude of the required network output while preserving the smooth interfacial transition structure.

We note that the proposed RDC is not specific to the Cahn--Hilliard equation. It can be applied more generally to partial differential equations whose solutions contain sharp transition layers between known asymptotic states.

\subsection{Variable scaling}
\label{subsec:variable_scaling_method}

When applying PINNs to MCL problems, the main training difficulty arises from the strongly localized and multiscale nature of the phase-field solution near the diffuse interface and the solid boundary. In particular, the order parameter $\phi$ exhibits sharp transitions across a thin interfacial layer, while the contact line region involves steep spatial gradients coupled with dynamic boundary conditions. These features worsen the ill-conditioning of the loss landscape, particularly when optimizing the Cahn--Hilliard residual and boundary residuals simultaneously. To address this challenge, we apply a \emph{variable-scaling} coordinate transformation to the governing equations prior to training \cite{ko2025vs}. The core idea is to rescale the coordinates associated with fine-scale interfacial and contact-line structures so that these localized variations are distributed over a larger computational range, thereby improving the stability and accuracy of PINNs training.

The method starts with a linear scaling of the independent variables. For the original computational domain in space and time, denoted as $(x, y, t) \in \Omega \times [t_{\text{min}}, t_{\text{max}}]$ (i.e., $t_{\text{min}}=t^n$, $t_{\text{max}}=t^{n+1}$ when applied within the time-marching scheme introduced in Section \ref{subsec:time_marching_scheme}, we define a set of scaled coordinates $(x', y', t')$ as:
\begin{equation} \label{eq:scaling_coeff_general}
    x' = N_x  x, \quad y' = N_y y, \quad t' = N_t t,
\end{equation}
where $N_x, N_y, N_t > 0$ are positive, tunable scaling hyperparameters. This transformation maps the original domain $\Omega \times [t_{\text{min}}, t_{\text{max}}]$ to a new scaled domain $\Omega' \times [N_t t_{\text{min}}, N_t t_{\text{max}}]$, where $\Omega'$ is the scaled spatial domain defined by the original bounds multiplied by the corresponding scaling factors.
Under this scaling, the differential operators transform via the chain rule as $\partial_x = N_x\,\partial_{x'}$, $\partial_x^2 = N_x^2\,\partial_{x'}^2$, with analogous expressions for $y$ and $t$.
Substituting into the governing equations \eqref{eq:CH_PDE1}, the initial condition \eqref{eq:CH_initial}, and the boundary conditions \eqref{eq:relaxation_bc}-\eqref{eq:conservation_bc}, the scaled problem becomes
\begin{equation} \label{eq:scaled_general_system}
\left\{
\begin{aligned}
    & N_t \frac{\partial \phi}{\partial t'} = M \left( N_x^2 \frac{\partial^2 \mu}{\partial x'^2} + N_y^2 \frac{\partial^2 \mu}{\partial y'^2} \right), && (x', y') \in \Omega', \\
    & \mu = -\varepsilon^2 \left( N_x^2 \frac{\partial^2 \phi}{\partial x'^2} + N_y^2 \frac{\partial^2 \phi}{\partial y'^2} \right) + F'(\phi), && (x', y') \in \Omega', \\
    & \phi\left(\frac{x'}{N_x}, \frac{y'}{N_y}, 0\right) = \phi^0\left(\frac{x'}{N_x}, \frac{y'}{N_y}\right), && (x', y') \in \Omega', \\
    & N_x \frac{\partial \phi}{\partial x'} n_{x'} + N_y \frac{\partial \phi}{\partial y'} n_{y'} = \frac{1}{\varepsilon} \left( \frac{\sqrt{2}\pi}{6} \cos(\theta_s) \cos\!\left(\frac{\pi}{2}\phi\right) - \frac{N_t}{\alpha} \frac{\partial \phi}{\partial t'} \right), && (x', y') \in \partial \Omega', \\
    & N_x \frac{\partial \mu}{\partial x'} n_{x'} + N_y \frac{\partial \mu}{\partial y'} n_{y'} = 0, && (x', y') \in \partial \Omega',
\end{aligned}
\right.
\end{equation}
where $(n_{x'},n_{y'})$ is the outward unit normal vector on $\partial \Omega'$.
The neural network is then trained to solve this modified system in the scaled coordinate space $(x', y', t')$. The scaling hyperparameters $N_x, N_y, N_t$ are problem-dependent and selected empirically. In this paper, $N_x$ and $N_y$ used in each numerical experiment are reported in Section~\ref{sec:experimental_results}, while $N_t$ is fixed to $1$. Once the training is completed, the solution predicted in the scaled coordinate space is mapped back to the original physical domain via the inverse transformation $x = x'/N_x$, $y = y'/N_y$, and $t = t'/N_t$ to recover the true physical solution.

\subsection{Adaptive loss weighting}
\label{subsec:Adaptive_loss_weighting}
In discrete-time PINNs, the composite loss \eqref{eq:ch_total_loss} consists of several residual terms associated with the governing equations and boundary conditions. The relative magnitudes of weights $\lambda_{r_\phi}$, $\lambda_{r_\mu}$, $\lambda_{b_\phi}$, and $\lambda_{b_\mu}$ play an important role in the training process. If the gradient contribution from one loss component is significantly larger than the others, the optimization may become biased toward that component, leading to an imbalanced enforcement of the governing equations and boundary conditions.

To reduce the dependence on manually selected loss weights, we adopt the uncertainty-based multi-task weighting strategy proposed by Kendall et al.~\cite{kendall2018multi} and subsequently adapted to PINNs by Hou et al.~\cite{hou2023enhancing}. In this approach, each loss component is modeled through a Gaussian likelihood with a learnable variance $\sigma_m^2$. Maximizing the corresponding joint log-likelihood leads to adaptive weights of the form $\lambda_m=e^{-s_m}$, where $s_m=\log\sigma_m^2$ and $m\in\{r_\phi,\,r_\mu,\,b_\phi,\,b_\mu\}$. For the four loss components considered here, the resulting objective is given by
\begin{equation} \label{eq:ch_weighted_loss}
    \mathcal{L}(\boldsymbol{\theta},\{s_m\})
    =
    \sum_{m\in\{r_\phi,\,r_\mu,\,b_\phi,\,b_\mu\}}
    \left(
        e^{-s_m}\mathcal{L}_m(\boldsymbol{\theta}) + s_m
    \right),
\end{equation}
where the log-variance parameters $s_m$ are trained together with the neural-network parameters $\boldsymbol{\theta}$. The factor $e^{-s_m}$ acts as an adaptive weight for $\mathcal{L}_m(\boldsymbol{\theta})$, whereas the regularization term $s_m$ penalizes excessive variance growth and prevents the corresponding loss component from being trivially down-weighted. This formulation enables an automatic balance between the governing-equation residuals and the boundary-condition residuals during training.

The composite loss function in \eqref{eq:ch_total_loss} consists of multiple competing objectives, whose gradient directions may be mutually conflicting. To improve the stability of the network optimization, we employ the Second-Order Aligned with Preconditioner (SOAP) optimizer \cite{vyas2024soap} to update the neural-network parameters $\boldsymbol{\theta}$. The adaptive weighting parameters $\{s_m\}$ are optimized separately using the Adam optimizer \cite{kingma2014adam} with an independent learning rate.
At each training iteration, the two parameter groups are updated in an alternating manner. Specifically, SOAP is first applied to update $\boldsymbol{\theta}$ by minimizing $\mathcal{L}(\boldsymbol{\theta})$, while the adaptive weights $\{s_m\}$ are kept fixed. Subsequently, Adam is used to update $\{s_m\}$ according to the likelihood-based objective $\mathcal{L}(\boldsymbol{\theta},\{s_m\})$, with $\boldsymbol{\theta}$ held fixed. This alternating optimization strategy separates the training of the neural-network parameters from the adaptation of the loss weights, thereby reducing unstable coupling between the two optimization processes.

\subsection{Adaptive sampling of collocation points} \label{subsec:Adaptive_sampling}

In the Cahn--Hilliard equation, the dynamics are dominated by thin interfacial layers where the order parameter transitions between the bulk values $\phi \approx \pm 1$. A uniform distribution of collocation points is therefore inefficient, as it allocates most resources to regions of nearly constant solution. To concentrate computational effort on the interfacial region, we employ an adaptive sampling strategy that selects collocation points based on the initial condition of each time interval.

The strategy is implemented by first identifying the diffuse interface region from the initial condition at $t=0$. Specifically, the initial phase field $\phi^0(x,y)$ is evaluated on a fine uniform grid, and the points away from the bulk phases are selected according to
\begin{equation} \label{eq:criterion_col_pts}
    |\phi^0(x,y)-1|>\tau
    \quad \text{and} \quad
    |\phi^0(x,y)+1|>\tau,
\end{equation}
where $\tau$ is a prescribed tolerance, taken as $\tau=0.1$ in this study.
This adaptive sampling procedure is applied at the beginning of each time interval in the time-marching scheme. For the network trained on $[t^n,t^{n+1}]$, the diffuse-interface region is re-extracted using the predicted solution at $t^n$ from the preceding network as the updated initial condition. In this way, the collocation points are refreshed at every time interval, allowing the sampling distribution to track the evolving interface.

\begin{figure}
\centering
\begin{subfigure}[b]{0.47\linewidth}
    \centering
    \includegraphics[width=\linewidth]{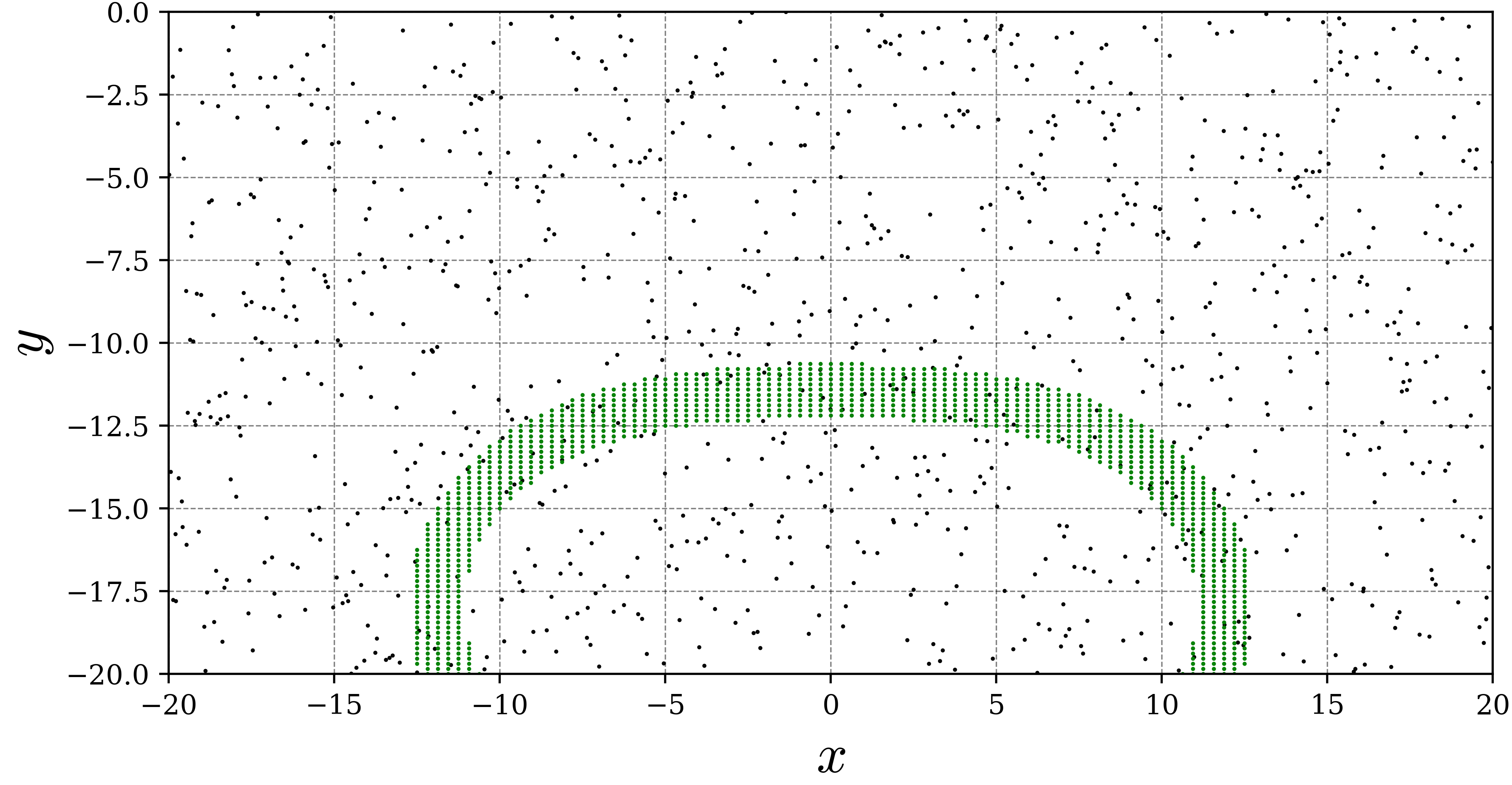}
    \caption{$t \in [0.9,\,1]$}
    \label{fig:AM_col_pts_1}
\end{subfigure}
\hfill
\begin{subfigure}[b]{0.47\linewidth}
    \centering
    \includegraphics[width=\linewidth]{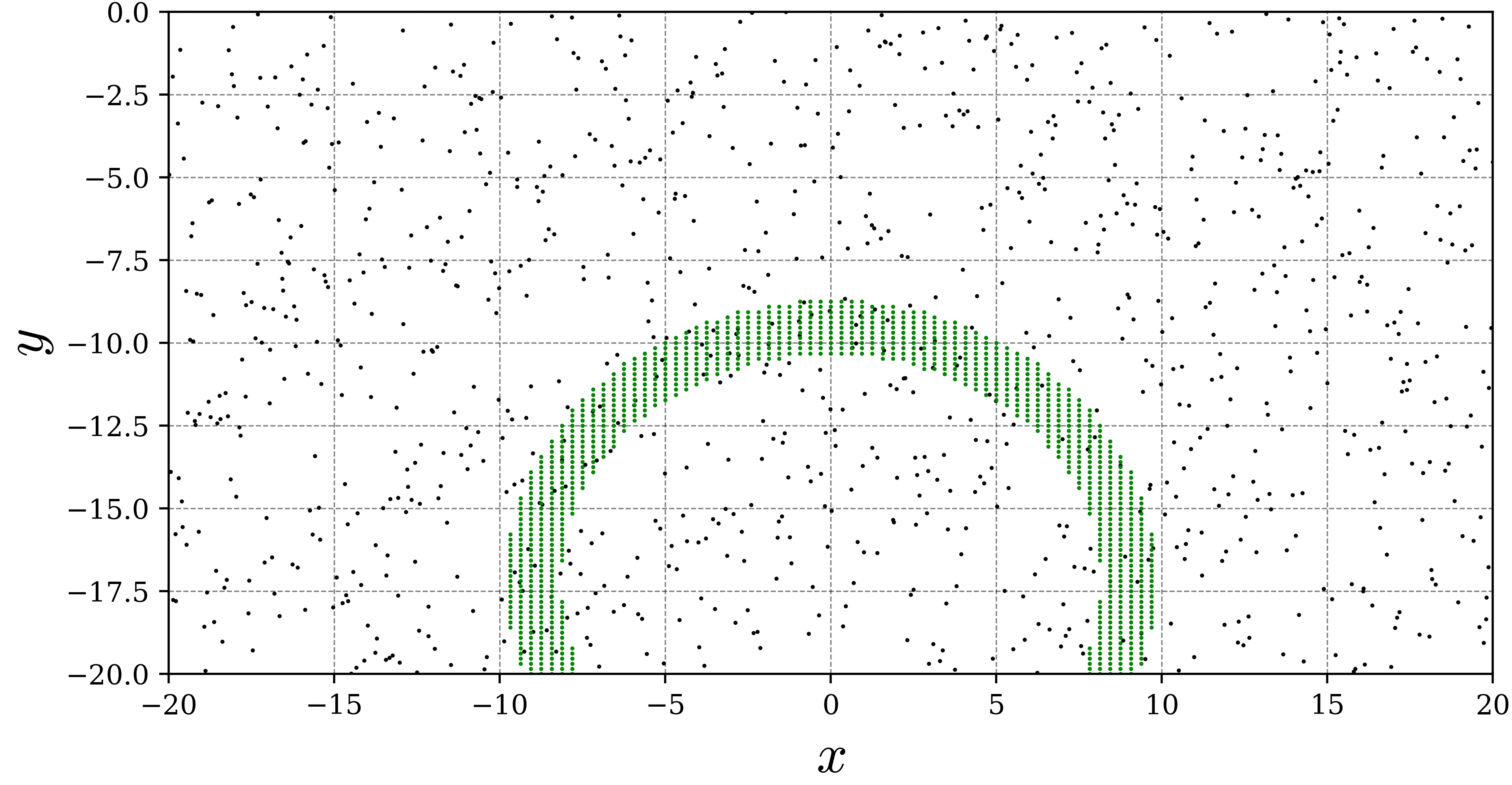}
    \caption{$t \in [3.9,\,4]$}
    \label{fig:AM_col_pts_4}
\end{subfigure}
\caption{Distribution of collocation points generated by the adaptive sampling method, illustrated using the droplet coalescence case at two selected time intervals. Black dots represent background points generated by Latin Hypercube Sampling (LHS), while green dots denote adaptive sampling points concentrated along the diffuse interface.}
\label{fig:adaptive_sampling}
\end{figure}

Figure~\ref{fig:adaptive_sampling} illustrates the spatial distribution of collocation points at two distinct time intervals. The collocation points consist of two different sets: the background points generated via Latin Hypercube Sampling (LHS), shown in black, and the adaptive points generated by the proposed sampling strategy, shown in green. The black points are sparsely distributed across the entire domain, whereas the green points are concentrated densely along the diffuse interface. As observed in the transition from $t \in [0.9, 1]$ to $t \in [3.9, 4]$, the cluster of green points dynamically shifts to track the propagation of the interface, demonstrating that the adaptive strategy effectively concentrates the computational efforts near the evolving interfacial region while the LHS points provide a sparse baseline coverage. 

\begin{remark}
The tolerance $\tau$ should be chosen in relation to the interface thickness parameter $\varepsilon$: a value $\tau \ll 1$ is needed to capture the full diffuse-interface profile, while excessively small $\tau$ may include near-bulk points and reduce sampling efficiency. In this work, $\tau = 10^{-1}$ provides a robust selection across all test cases. We note that the sampling band defined by this criterion is much wider than the diffuse interface itself.
Since the collocation points are determined from the initial condition of each time interval and remain fixed throughout that interval, the wide sampling band provides a spatial buffer that accommodates interface displacement, and the continuous parameterization of the neural network enables it to generalize to nearby regions even where collocation points are sparser. The high-order IRK discretization ($q \gg 1$) further provides fine temporal resolution within each interval, allowing the network to track the evolving interface through its stage values.
\end{remark}

\subsection{Symmetry preservation through neural network inputs}
\label{subsec:symmetry_embedding}

When the governing equations, initial data, and boundary conditions are invariant under a spatial transformation, the exact solution is expected to preserve the same symmetry. To incorporate this property, we enforce \emph{symmetry preservation} at the level of the neural network inputs by constructing input transformations that are invariant under the prescribed spatial symmetry. This architectural design embeds the symmetry prior directly into the network representation and helps prevent the optimizer from exploring nonphysical, symmetry-breaking configurations.

In the droplet coalescence problem considered in Section~\ref{sec:experimental_results}, the initial condition (two nuclei symmetric with respect to the $y$-axis) and the boundary conditions are both invariant under the reflection $x \mapsto -x$, so that $\phi(-x, y, t) = \phi(x, y, t)$ for all $t \geq 0$. To embed this even symmetry, we preprocess the $x$-coordinate via the mapping
\begin{equation}
\label{eq:input_mapping}
x \mapsto N_x \cos\!\left( \frac{x}{N_x} \right),
\end{equation}
where $N_x$ is the scaling factor along the $x$-direction used in the variable-scaling method. Since $\cos(-z) = \cos(z)$, the network receives identical inputs for any pair $(x, y)$ and $(-x, y)$, so that $\phi(x, y, t) = \phi(-x, y, t)$ by construction. The $y$-coordinate is passed directly to the network without transformation.

\begin{remark}
This symmetry-preserving input transformation is applied only when the computational domain and the solution are both invariant under the reflection $x \mapsto -x$. In particular, the physical domain $[x_{\min},x_{\max}]$ must be symmetric with respect to $x=0$, i.e., $x_{\min}=-x_{\max}$, so that the problem can be reduced to the half-domain $[0,x_{\max}]$. On this reduced domain, the mapping $x \mapsto N_x\cos(x/N_x)$ is one-to-one provided that $[0,x_{\max}] \subset [0,\pi N_x]$. If either the computational domain is not symmetric about $x=0$, or the solution does not possess even symmetry in $x$, such as the shear-induced deformation problem in Section~\ref{sec:experimental_results}, this input transformation is not applied.
\end{remark}

Algorithm~\ref{alg:mcl_pinns} summarizes the proposed neural solver MCL-PINNs, integrating all components described in the preceding subsections.

\begin{algorithm}
\caption{A phase-field neural solver for moving contact line problems (MCL-PINNs)}\label{alg:mcl_pinns}
\begin{algorithmic}[1]
\REQUIRE Initial condition $\phi^0(x,y)$; final time $T$; number of time intervals $N$; IRK coefficients $\{a_{ij}, b_j, c_j\}$; scaling parameters $N_x, N_y, N_t$; tolerance $\tau$; number of training iterations $N_{\text{iter}}$.
\ENSURE Sequence of trained networks $\{\mathcal{N}_n\}_{n=0}^{N-1}$ approximating $\phi$ on $\Omega \times [0,T]$.
\STATE  Set $\Delta t \leftarrow T/N$.
\STATE \emph{Variable scaling}: scale the system using $N_x, N_y, N_t$ via \eqref{eq:scaling_coeff_general}.
\STATE \emph{Time-marching scheme}:
\FOR{$n = 0, 1, \ldots, N-1$}
    \STATE \emph{Adaptive sampling:} Identify interfacial collocation points from $\phi^n$ using criterion~\eqref{eq:criterion_col_pts} with tolerance $\tau$; combine with background points.
    \STATE Initialize network $\mathcal{N}_n$: if $n = 0$, initialize $\boldsymbol{\theta}$ randomly; otherwise, warm-start by inheriting the parameters of $\mathcal{N}_{n-1}$. Apply the \emph{relaxed distribution constraint}~\eqref{eq:RDC} to the outputs, if applicable. 
    Apply the \emph{symmetry preservation} mapping~\eqref{eq:input_mapping} to the inputs, if applicable.
    \STATE Initialize \emph{adaptive loss weighting} parameters $\{s_m\}$, $m \in \{r_\phi, r_\mu, b_\phi, b_\mu\}$.
    \FOR{$\text{iter} = 1, \ldots, N_{\text{iter}}$}
    \STATE Evaluate the loss components 
    $\mathcal{L}_{r_{\phi}}$, $\mathcal{L}_{r_{\mu}}$, 
    $\mathcal{L}_{b_{\phi}}$, and $\mathcal{L}_{b_{\mu}}$.
    \STATE Form the weighted network loss 
    $\mathcal{L}(\boldsymbol{\theta})$
    using the current adaptive weights~\eqref{eq:ch_total_loss}.
    \STATE Take one SOAP step to update $\boldsymbol{\theta}$ using 
    $\nabla_{\boldsymbol{\theta}}\mathcal{L}(\boldsymbol{\theta})$, 
    while keeping $\{s_m\}$ fixed.
    \STATE Form the likelihood-based weighting objective 
    $\mathcal{L}(\boldsymbol{\theta};\{s_m\})$ ~\eqref{eq:ch_weighted_loss}.
    \STATE Take one Adam step to update $\{s_m\}$ using 
    $\nabla_{\{s_m\}}\mathcal{L}(\boldsymbol{\theta}; \{s_m\})$, 
    while keeping $\boldsymbol{\theta}$ fixed.
    \ENDFOR
\STATE Evaluate the terminal prediction $\tilde{\phi}^{n+1}$ from the trained network $\mathcal{N}_n$ and set 
$\phi^{n+1} \leftarrow \tilde{\phi}^{n+1}$.
\ENDFOR
\end{algorithmic}
\end{algorithm}

\section{Experimental results} \label{sec:experimental_results}
In this section, we assess the performance of MCL-PINNs through three representative problems, namely droplet coalescence, shear-induced droplet deformation, and dynamic wetting through a heterogeneous channel. The MCL-PINNs predictions are quantitatively evaluated against finite element method (FEM) reference solutions using two standard metrics: the mean squared error (MSE) and the relative $L_2$ error.
These metrics are defined as
\begin{equation}
    \text{MSE} = \frac{1}{N_{test}} \sum_{i=1}^{N_{test}} \left| \phi_{\text{MCL-PINNs}}(\mathbf{x}_i) - \phi_{\text{FEM}}(\mathbf{x}_i) \right|^2,
    \label{eq:mse}
\end{equation}
\begin{equation}
    \text{Relative } L_2 \text{ Error} = \frac{\sqrt{\sum_{i=1}^{N_{test}} \left| \phi_{\text{MCL-PINNs}}(\mathbf{x}_i) - \phi_{\text{FEM}}(\mathbf{x}_i) \right|^2}}{\sqrt{\sum_{i=1}^{N_{test}} \left| \phi_{\text{FEM}}(\mathbf{x}_i) \right|^2}},
    \label{eq:rel_l2}
\end{equation}
where $N_{test}$ denotes the total number of test points, $\mathbf{x}_i$ represents the spatial coordinates of the $i$-th evaluation point, and $\phi_{\text{MCL-PINNs}}$ and $\phi_{\text{FEM}}$ denote the phase field solutions predicted by MCL-PINNs and computed by the FEM reference solver, respectively. For the FEM reference solver, we use continuous piecewise linear ($P_1$) finite elements on a structured mesh, with the mesh size set to one-half of the interface thickness to adequately resolve the interfacial evolution. An energy-stable second-order accurate scheme \cite{zhang2025fully} is employed for time integration.

For all numerical experiments presented in this work, a consistent implementation framework is employed to ensure fair comparison and reproducibility. The MCL-PINNs architecture comprises 6 fully connected hidden layers, each containing 128 neurons. A fixed time-marching step $\Delta t = 0.1$ is applied for all simulations. A fixed learning rate of $10^{-3}$ is used for the SOAP optimizer to update the neural network parameters $\boldsymbol{\theta}$. For the adaptive loss weighting mechanism, the learnable log-variance parameters $\{s_m\}_{m\in\{r_\phi,\,r_\mu,\,b_\phi,\,b_\mu\}}$ are initialized as $s_{r_\phi}=-\log(1000)$ and $s_{b_\phi}=s_{r_\mu}=s_{b_\mu}=0$, and are 
optimized using the Adam optimizer with a fixed learning rate of $10^{-3}$. The selected interface-focused collocation points are chosen according to criterion~\eqref{eq:criterion_col_pts} with $\tau=0.1$. All numerical experiments are implemented using PyTorch 2.0.0 and executed on a workstation equipped with a single NVIDIA GeForce RTX 4090 GPU.

\subsection{Droplet coalescence}
%------------------------------------------------------------------
We first examine the droplet coalescence problem to validate the effectiveness and accuracy of MCL-PINNs. The computational domain is $\Omega = [-1, 1] \times [-1, 0]$ and $t \in [0, 8]$. The initial setup consists of two circular droplets centered at $(\pm 0.28, -1)$ with a radius of $r = 0.4$. The initial condition is given by
$$
\phi^0(x, y) = \max\left\{ 
\tanh\left( \frac{0.4 - \sqrt{(x - 0.28)^2 + (y + 1)^2}}{2\varepsilon} \right), 
\tanh\left( \frac{0.4 - \sqrt{(x + 0.28)^2 + (y + 1)^2}}{2\varepsilon} \right)
\right\},
$$
where the interface thickness parameter is set to $\varepsilon = 0.02$. Dynamic boundary conditions \eqref{eq:relaxation_bc}-\eqref{eq:conservation_bc} are imposed on all boundaries.
Two representative static contact angles are considered: $\theta_s = 50^\circ$ and $\theta_s = 120^\circ$, representing hydrophilic and hydrophobic wetting regimes. The mobility coefficient and contact line relaxation parameter are set to $M = 1$ and $\alpha = 100$, respectively. The external force is neglected in this problem, i.e., $\mathbf{u}=\mathbf{0}$.

A $q$-stage IRK scheme with a fixed time step $\Delta t = 0.1$ is used for the temporal discretization. The number of stages is adjusted to address the numerical stiffness of each case: $q = 20$ for $\theta_s = 50^\circ$ and $q = 50$ for $\theta_s = 120^\circ$. The relaxed distribution constraint \eqref{eq:RDC} with $c=0.9$ is applied to the network output, and the symmetry preservation mapping \eqref{eq:input_mapping} is used for the input in this problem. A variable-scaling transformation is applied with scaling factors $N_x = 20$, $N_y = 20$, and $N_t = 1$. The interior collocation points are generated using a hybrid approach that combines interface-aware location with randomly distributed interior samples. For the initial time interval, a uniform candidate grid of $50 \times 50$ nodes is generated over the domain. In subsequent time intervals, this selection process is applied to a refined candidate grid of $129 \times 129$ nodes to better capture the evolving interface with higher spatial resolution. These interface-focused points are then combined with 1,000 additional interior points distributed via LHS, ensuring the capture of interface dynamics while maintaining sufficient coverage of the bulk phases. Moreover, we use 41 uniformly distributed nodes along the horizontal boundaries and 21 nodes along the vertical boundaries. The number of training iterations is $N_{\text{iter}}=10,000$ for case $50^\circ$ and $N_{\text{iter}}=20,000$ for case $120^\circ$, respectively. 

\begin{figure}
    \centering
    \begin{subfigure}[b]{0.48\textwidth}
        \centering
        \includegraphics[width=\textwidth]{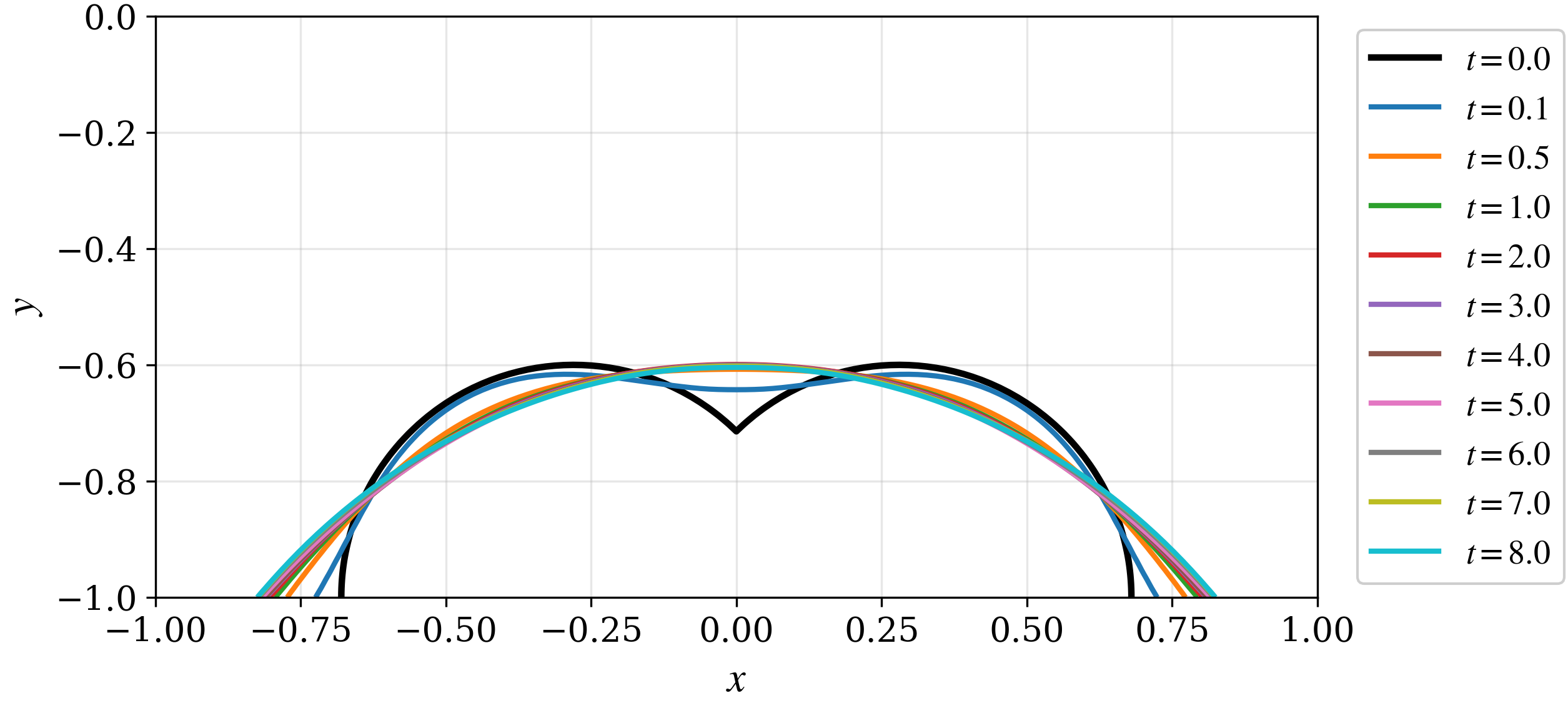}
        \caption{$\theta_s=50^\circ$}
        \label{fig:phi_0_evo_50}
    \end{subfigure}
    \hfill
    \begin{subfigure}[b]{0.48\textwidth}
        \centering
        \includegraphics[width=\textwidth]{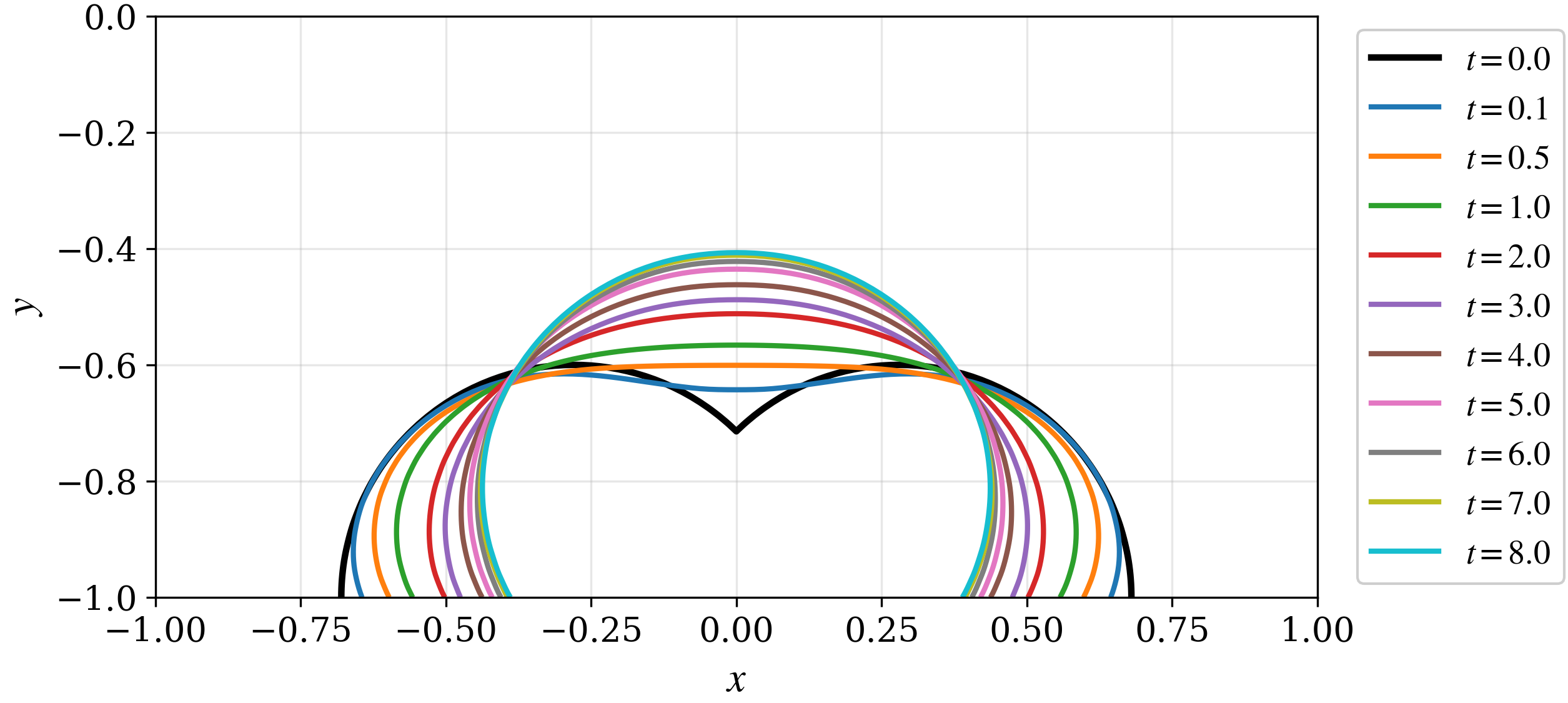}
        \caption{$\theta_s=120^\circ$}
        \label{fig:phi_0_evo_120}
    \end{subfigure}
    \caption{Evolution of the interface ($\phi=0$) using MCL-PINNs for the droplet coalescence problem with different static contact angles.}
    \label{fig:phi_evolution}
\end{figure}

Figure~\ref{fig:phi_evolution} illustrates the temporal evolution of the interface ($\phi=0$) for both contact angle configurations. For different static contact angles, the droplet moves away from its initial position in different directions under the effect of the unbalanced Young stress and the surface tension, and eventually reaches a stable state after attaining the prescribed static contact angle. For the hydrophilic case ($\theta_s = 50^\circ$), the coalesced droplet spreads more readily along the substrate, forming a flatter equilibrium profile with a smaller contact angle. In contrast, the hydrophobic case ($\theta_s = 120^\circ$) exhibits a more pronounced spherical cap shape with limited spreading. The coalescence process completes around $t \approx 3.0$-$4.0$, after which the merged droplet gradually relaxes toward its equilibrium configuration. The smooth evolution of the interface demonstrates the capability of MCL-PINNs to capture complex topological changes and moving contact line dynamics.
%------------------------------------------------------------------
%------------------------------------------------------------------

%------------------------------------------------------------------

\subsubsection{Ablation study} \label{subsubsec:case1_ablation}
To systematically evaluate the contribution of each component in the MCL-PINNs framework, we conduct comprehensive ablation experiments. The standard PINNs \cite{raissi2019physics} and a time-marching PINNs (TM-PINNs) \cite{wight2020solving} are also used for comparison. Table~\ref{tab:ablation} presents the quantitative comparison of MSE and relative $L_2$ errors for different method configurations across both contact angle cases. It can be observed that standard PINNs without time-marching result in $\mathcal{O}(1)$ mean squared and relative $L_2$ errors, confirming that temporal decomposition is essential for long-time integration. The time-marching scheme alone reduces errors by roughly one order of magnitude yet still leaves notable residuals near interfaces.

The complete MCL-PINNs framework achieves mean squared errors of $2.49\times10^{-4}$ for the $50^\circ$ case and $1.83\times10^{-3}$ for $120^\circ$, with corresponding relative $L_2$ errors of $1.63\times10^{-2}$ and $4.41\times10^{-2}$. These results represent significant improvement over baseline PINNs and TM-PINNs. Among individual components, removing adaptive sampling causes the largest performance drop, increasing mean squared error by over 80 times at $50^\circ$ and roughly 180 times at $120^\circ$, which highlights its critical role in resolving sharp interfaces. Ablating variable scaling or symmetry preservation raises errors by 20 to 30 times, while removing relaxed distribution constraint or disabling adaptive weighting leads to more moderate but still noticeable degradation.
We also notice that the hydrophobic configuration at $120^\circ$ consistently exhibits greater sensitivity to component removal. This suggests that hydrophobic wetting dynamics, characterized by more active contact line motion, impose higher demands on the stabilization of prediction.
% 120 (version num: 42.5.17.15.22)
\begin{table}
\centering
\caption{Ablation study comparing the performance of different method configurations for the droplet coalescence problem at $t = 8$. TM: time-marching scheme; AS: adaptive sampling strategy; VS: variable scaling scheme; SYM: symmetry preservation in neural network inputs; RDC: Relaxed distribution constraint; AW: adaptive weighting method. MSE: Mean squared error; Rel. err.: Relative $L_2$ error. The notation ``(w/o X)'' denotes the ablation variant without using technique ``X'', while the rest of the proposed framework remains unchanged.}
\label{tab:ablation}
\footnotesize
\setlength{\tabcolsep}{2pt}
\begin{tabularx}{\textwidth}{@{} l *{6}{C} *{4}{c} @{}}
\toprule
\multirow{2}{*}{Method} & \multirow{2}{*}{TM} & \multirow{2}{*}{AS} & \multirow{2}{*}{VS} & \multirow{2}{*}{SYM} & \multirow{2}{*}{RDC} & \multirow{2}{*}{AW} & \multicolumn{2}{c}{$\theta_s=50^\circ$} & \multicolumn{2}{c}{$\theta_s=120^\circ$} \\
\cmidrule(lr){8-9} \cmidrule(lr){10-11}
 & & & & & & & MSE & Rel. err. & MSE & Rel. err. \\
\midrule
PINNs & \xmark & \xmark & \xmark & \xmark & \xmark & \xmark & 1.14 & 1.10 & 1.77 & 1.37 \\
TM-PINNs & \cmark & \xmark & \xmark & \xmark & \xmark & \xmark & 8.97e-1 & 9.75e-1 & 3.85e-1 & 6.40e-1\\
MCL-PINNs & \cmark & \cmark & \cmark & \cmark & \cmark & \cmark & 2.49e-4 & 1.63e-2 & 1.83e-3 & 4.41e-2 \\
MCL-PINNs (w/o TM) & \xmark & \cmark & \cmark & \cmark & \cmark & \cmark & 9.29e-1 & 9.93e-1 & 8.33e-1 & 9.42e-1\\
MCL-PINNs (w/o AS) & \cmark & \xmark & \cmark & \cmark & \cmark & \cmark & 2.03e-2 & 1.47e-1 & 3.23e-1 & 5.86e-1\\
MCL-PINNs (w/o VS) & \cmark & \cmark & \xmark & \cmark & \cmark & \cmark & 6.83e-4 & 2.69e-2 & 5.48e-2 & 2.41e-1\\
MCL-PINNs (w/o SYM) & \cmark & \cmark & \cmark & \xmark & \cmark & \cmark & 9.66e-3 & 1.01e-1 & 3.62e-2 & 1.96e-1\\
MCL-PINNs (w/o RDC) & \cmark & \cmark & \cmark & \cmark & \xmark & \cmark & 6.41e-3 & 8.24e-2 & 2.05e-2 & 1.48e-1\\
MCL-PINNs (w/o AW) & \cmark & \cmark & \cmark & \cmark & \cmark & \xmark & 2.11e-3 & 4.73e-2 & 1.26e-1 & 3.67e-1\\
\bottomrule
\end{tabularx}
\end{table}

%------------------------------------------------------------------
\begin{figure}
    \centering
    \includegraphics[width=\textwidth,height=0.85\textheight,keepaspectratio]{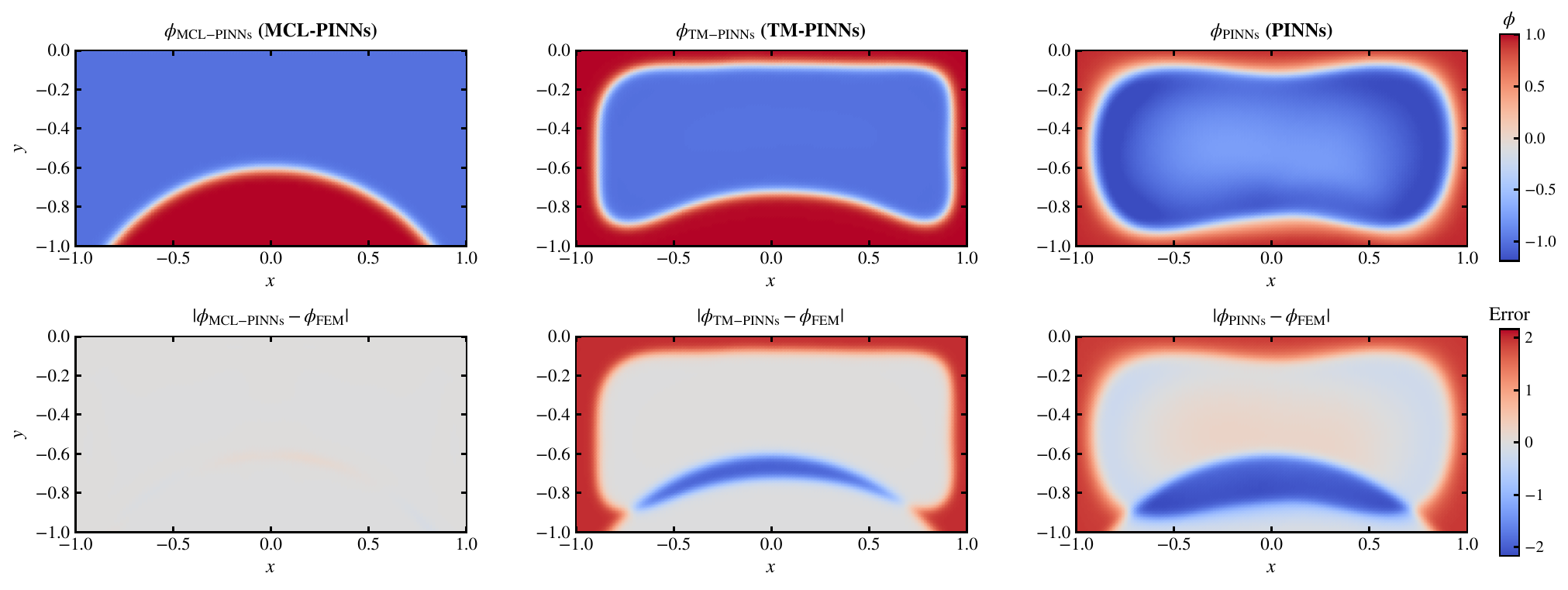}
    \caption{Comparison of results from MCL-PINNs, TM-PINNs, and standard PINNs at $t = 8$ for case $\theta_s = 50^\circ$ of the droplet coalescence problem. Top row: phase-field distributions. Bottom row: absolute error distributions.}
    \label{fig:comparison_50}
\end{figure}
%------------------------------------------------------------------
\begin{figure}
    \centering
    \includegraphics[width=\textwidth,height=0.85\textheight,keepaspectratio]{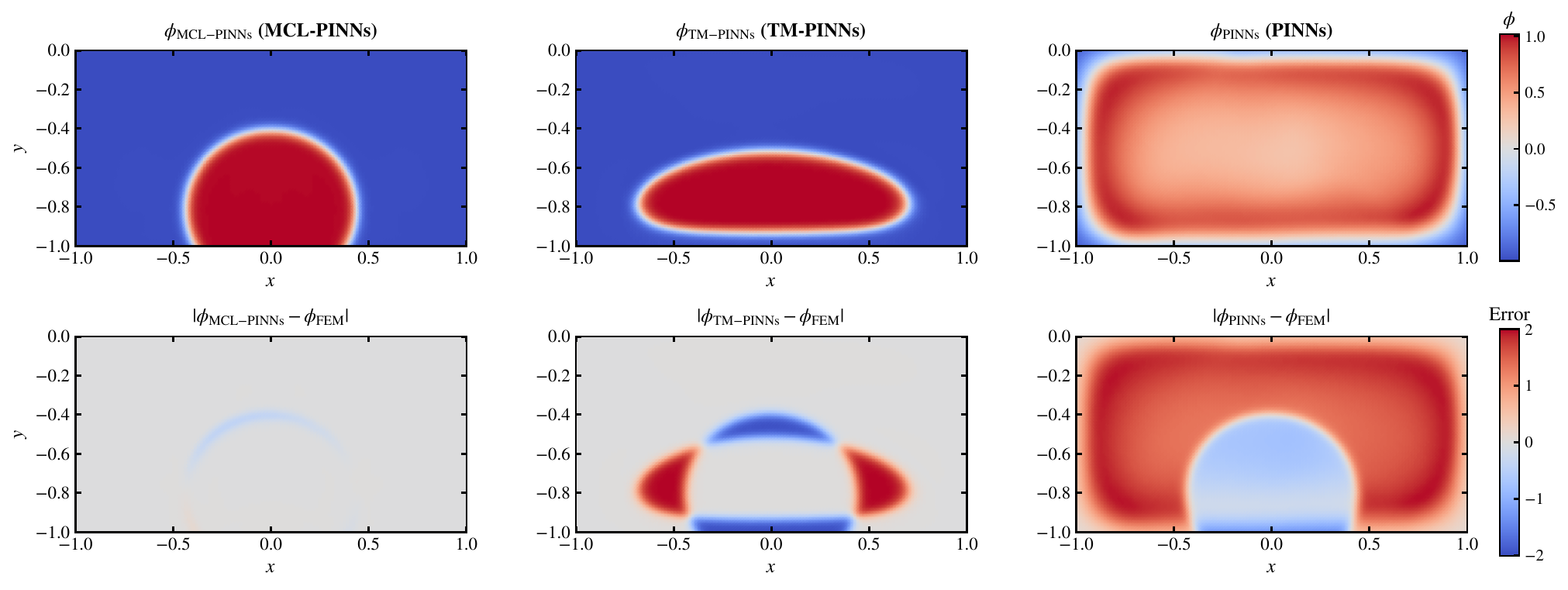}
    \caption{Comparison of results from MCL-PINNs, TM-PINNs, and standard PINNs at $t = 8$ for case $\theta_s = 120^\circ$ of the droplet coalescence problem. Top row: phase-field distributions. Bottom row: absolute error distributions.}
    \label{fig:comparison_120}
\end{figure}
%------------------------------------------------------------------

Figures~\ref{fig:comparison_50} and~\ref{fig:comparison_120} present comparisons of the results at the final time step $t = 8$, when the system reaches near-equilibrium. For the case $\theta_s=50^\circ$ (Figure~\ref{fig:comparison_50}), MCL-PINNs produces phase-field distributions almost identical to the FEM reference solution. The error distribution reveals that deviations are primarily confined to the immediate interfacial region. In contrast, TM-PINNs exhibits noticeable errors near the  interface and the boundaries, while standard PINNs fails to capture the correct equilibrium configuration, producing a completely incorrect droplet shape with substantial phase mixing.
The case $\theta_s=120^\circ$ in Figure~\ref{fig:comparison_120} presents more demanding dynamics due to the pronounced curvature. MCL-PINNs successfully captures the spherical cap morphology, demonstrating robustness across different wetting regimes. 
TM-PINNs shows significant deviation in the predicted contact line position, whereas standard PINNs produces a severely distorted interface profile.

%------------------------------------------------------------------
\begin{figure}
    \centering
    \includegraphics[width=1\linewidth]{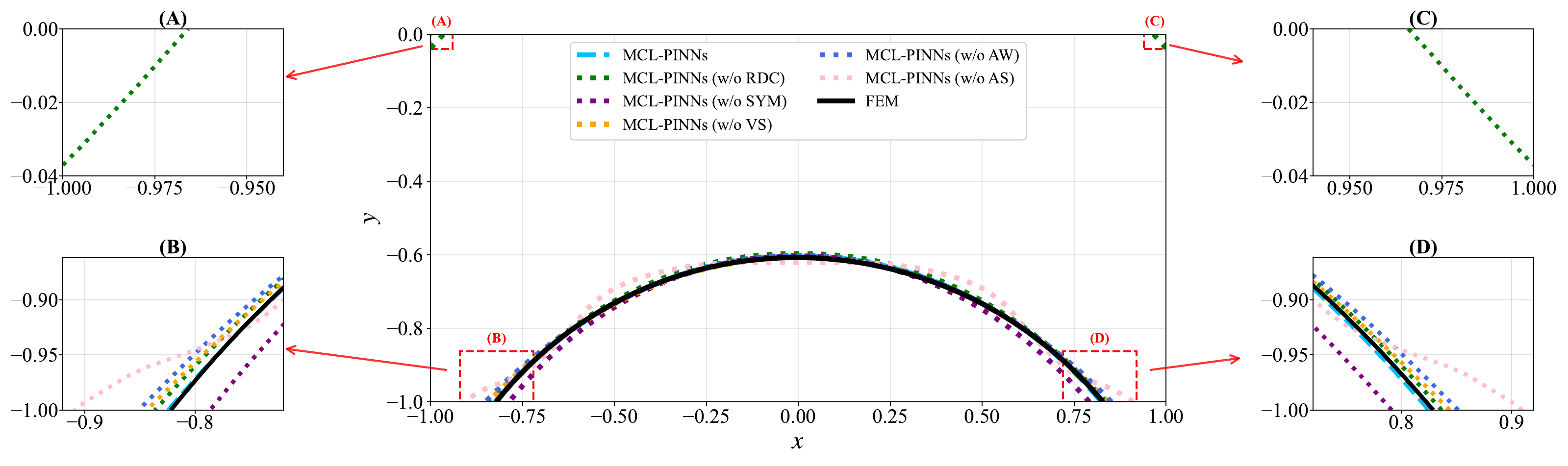}
    \caption{A comparison of the interface ($\phi=0$) at $t=8$ for case $\theta_s=50^\circ$ of the droplet coalescence problem. Panels (A) and (C) provide enlarged views of the upper-left and upper-right corner regions, respectively. Panels (B) and (D) show enlarged views of the contact line region near the solid substrate.}
    \label{fig:interface_zoom_50}
\end{figure}
%------------------------------------------------------------------
\begin{figure}
    \centering
    \includegraphics[width=1\linewidth]{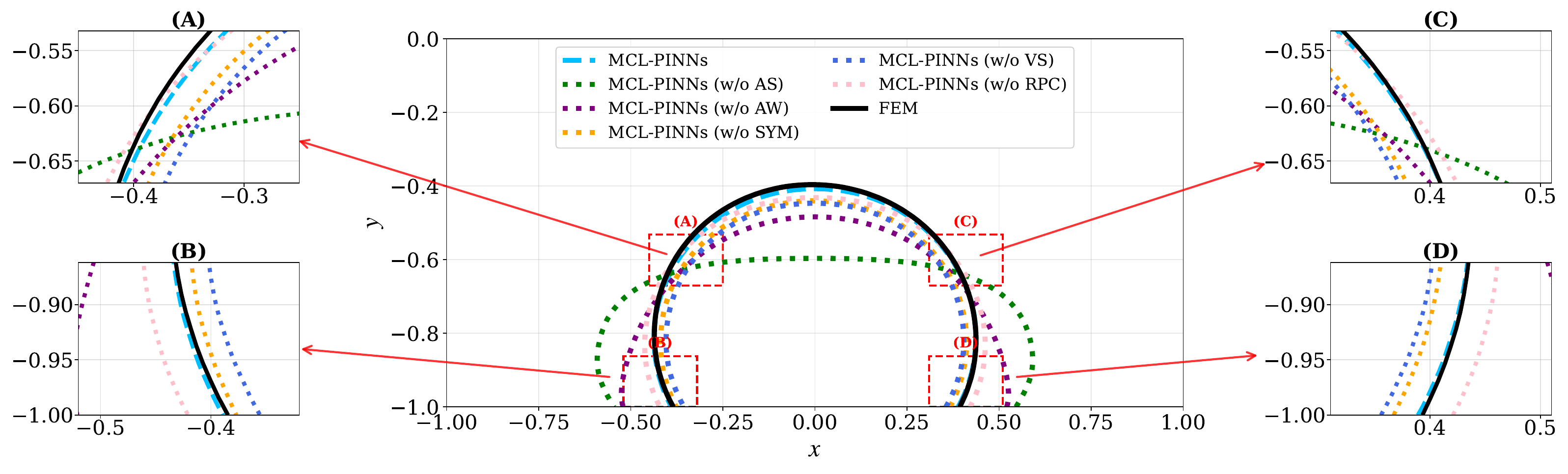}
    \caption{A comparison of the interface ($\phi=0$) at $t=8$ for case $\theta_s=120^\circ$ of the droplet coalescence problem.  Panels (A) and (C) enlarge the views near the middle part of the interface, while panels (B) and (D) provide enlarged views near the solid boundary.}
    \label{fig:interface_zoom_120}
\end{figure}
%------------------------------------------------------------------
The detailed interface comparisons in Figure~\ref{fig:interface_zoom_50} and Figure~\ref{fig:interface_zoom_120} serve as direct visual evidence of the predictive accuracy and structural robustness of the proposed MCL-PINNs framework. The complete MCL-PINNs framework demonstrates good agreement with the FEM reference solution, with the interfaces nearly indistinguishable across all regions. In contrast, the ablated variants exhibit varying degrees of deviation from the benchmark, confirming that each integrated component is necessary for maintaining quantitative accuracy. Specifically, for the hydrophilic case ($\theta_s = 50^\circ$, Figure~\ref{fig:interface_zoom_50}), the enlarged views in panels (A) and (C)  demonstrate a critical failure when the RDC is not used, which generates spurious interface artifacts in these regions. This phenomenon arises because, in the absence of RDC, the neural network is free to produce unphysical overshoots and undershoots, violating the bounds $\phi \in [-1, 1]$. 
For the hydrophobic case ($\theta_s = 120^\circ$, Figure~\ref{fig:interface_zoom_120}), the enlarged views in panels (A)--(D) show that the ablated variants lead to great deviations. In particular, the configurations without adaptive weighting or adaptive sampling of collocation points lead to substantial shifts in the contact position and interface profile, indicating that these two components are crucial for resolving hydrophobic cases.
These comparisons demonstrate that the combined use of all proposed components establishes MCL-PINNs as a robust and reliable framework for simulating moving contact line dynamics.
%------------------------------------------------------------------

\subsubsection{Effect of the relaxed distribution constraint factor}
\label{subsub:rdc_sensitivity}
\begin{figure}
\centering
\includegraphics[width=1\textwidth]{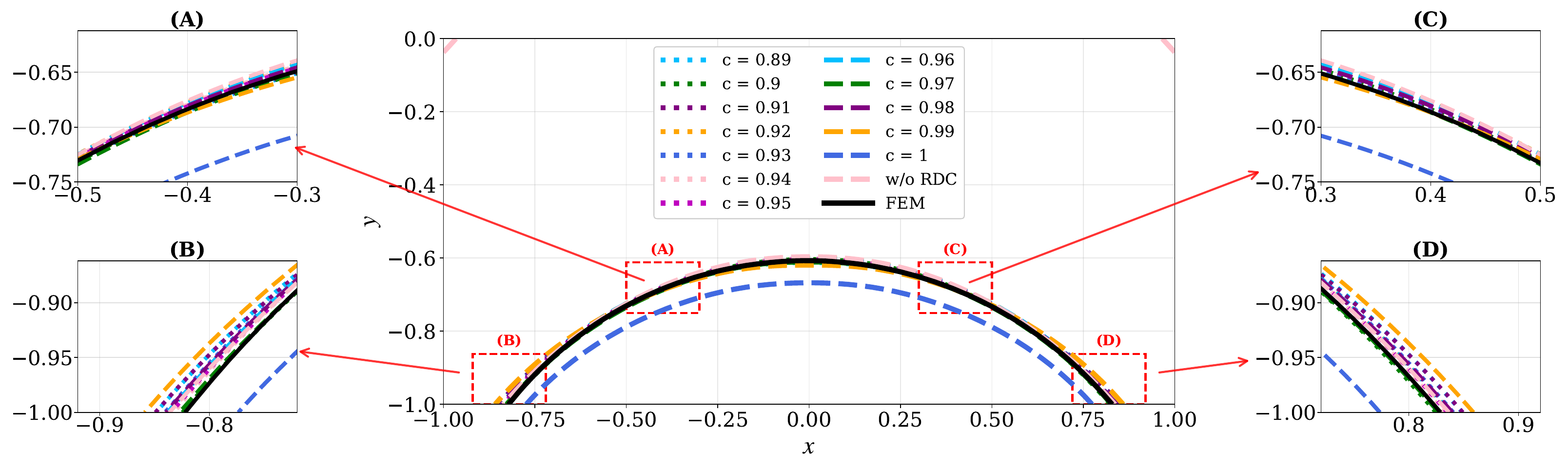}
\caption{A comparison of the interface ($\phi=0$) obtained by using the relaxed distribution constraint (RDC) method in MCL-PINNs with different values of the factor $c$, MCL-PINNs without RDC (`w/o RDC'), and the FEM solution at $t=8$ for case $\theta_s=50^\circ$ of the droplet coalescence problem. Panels (A)--(D) present zoomed-in views of selected regions.}
\label{fig:tanh_c_case1}
\end{figure}

To investigate the influence of relaxed distribution constraint (RDC) on the performance of the neural solver, we conduct a sensitivity analysis with respect to the scaling factor $c$ in \eqref{eq:RDC} for the droplet coalescence problem with $\theta_s=50^\circ$. Figure~\ref{fig:tanh_c_case1} shows a comparison of the interface ($\phi=0$) at $t=8$ across different values of $c$. Table~\ref{tbl:c_sensitivity_row} reports the MSE and relative $L_2$ error between the neural network predictions and the FEM reference solutions. As observed, although the standard profile constraint with $c=1.0$ suppresses the spurious interfacial artifacts near the top corners, it yields a larger error than the configuration without this technique. This suggests that enforcing $\tilde{\phi}=\tanh(z)$ without relaxation introduces an overly restrictive output representation, since the bulk values $\tilde{\phi}\approx\pm 1$ require large pre-activation values. 
In contrast, introducing the relaxation factor $c<1$ significantly mitigates this issue by confining the admissible output range of the network to $(-\operatorname{arctanh}(c), \operatorname{arctanh}(c))$. These results indicate that choosing $c$ in the range $0.90$--$0.97$ is effective for resolving sharp interfacial dynamics, with relatively low errors observed at $c=0.90$ and $c=0.97$. Since a smaller value of $c$ further reduces the magnitude of the required pre-activation values, it provides a more favorable output representation for network training. Therefore, we set $c=0.90$ in the subsequent experiments.

\begin{table}[htbp]
\centering
\caption{Computational error at $t = 8$ between the neural network solution and FEM reference under varying relaxation factors $c$ in the droplet coalescence case with $\theta_s=50^\circ$. MSE: Mean squared error; Rel. err.: Relative $L_2$ error.}
\label{tbl:c_sensitivity_row}
\footnotesize
\resizebox{\textwidth}{!}{%
\begin{tabular}{l *{12}{c}} 
\toprule
$c$ & 0.89 & 0.90 & 0.91 & 0.92 & 0.93 & 0.94 & 0.95 & 0.96 & 0.97 & 0.98 & 0.99 & 1.00 \\
\midrule
MSE       & 2.32e-3 & 2.49e-4 & 2.11e-3 & 2.51e-3 & 1.27e-3 & 1.27e-3 & 8.25e-4 & 2.16e-3 & 1.51e-4 & 1.34e-3 & 6.20e-3 & 9.15e-2 \\
Rel. err. & 4.96e-2 & 1.63e-2 & 4.73e-2 & 5.16e-2 & 3.67e-2 & 3.66e-2 & 2.96e-2 & 4.79e-2 & 1.26e-3 & 3.77e-2 & 8.11e-2 & 3.12e-1 \\
\bottomrule
\end{tabular}
}
\end{table}

\subsubsection{Temporal accuracy and stability}
\label{subsubsec:case1_stability}
%------------------------------------------------------------------
\begin{figure}
    \centering
    \includegraphics[width=1\textwidth]{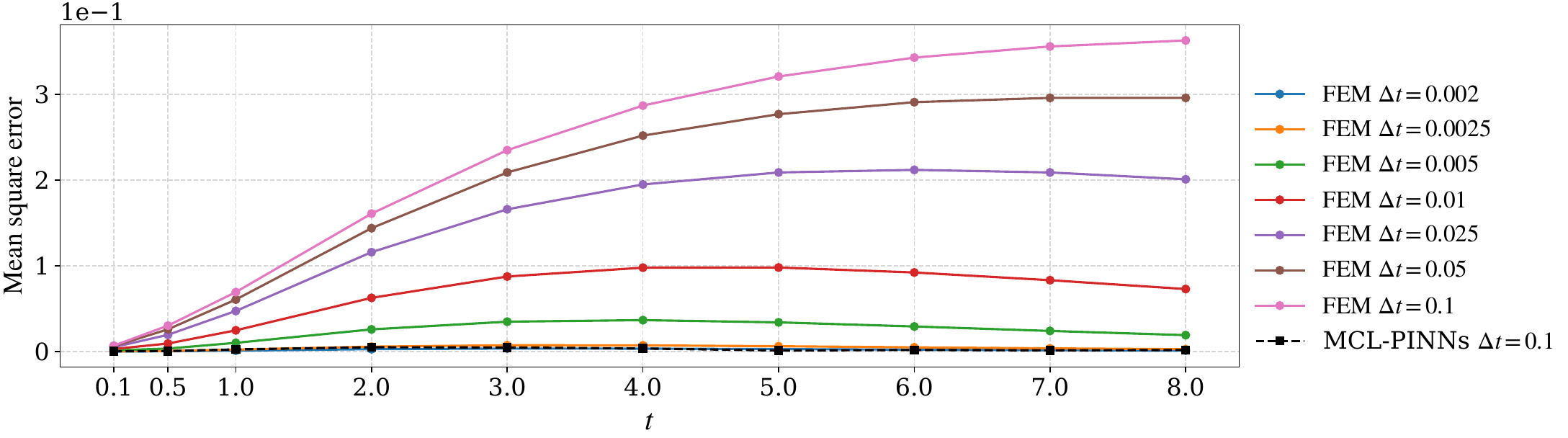}
    \caption{Temporal evolution of mean squared error (MSE) for case $\theta_s = 120^\circ$ of the droplet coalescence problem. The proposed MCL-PINNs framework with $\Delta t = 0.1$ achieves accuracy comparable to the reference FEM solution with $\Delta t = 0.001$, while FEM solutions with coarser time steps exhibit pronounced error accumulation.}
    \label{fig:mse_time_120du_FEM_vs_PINNs}
\end{figure}

Figure~\ref{fig:mse_time_120du_FEM_vs_PINNs} provides a detailed comparison of the temporal evolution of mean squared error (MSE) between the MCL-PINNs framework and the FEM solver for the hydrophobic case ($\theta_s = 120^\circ$). All errors are measured relative to the reference FEM solution obtained with the finest time step $\Delta t = 0.001$, which serves as the accuracy benchmark.
A notable distinction emerges in the error propagation of the two approaches. FEM solutions display a clear accumulation of error over time, with MSE values increasing monotonically throughout the simulation. This effect is especially pronounced for larger time step sizes $\Delta t \geq 0.025$, where the error grows nearly an order of magnitude from $t=0$ to $t=8$. Even for small time step sizes $\Delta t = 0.002$--$0.01$, the cumulative error exhibits a consistent upward trend. In contrast, the MCL-PINNs framework with $\Delta t = 0.1$ demonstrates high temporal accuracy and stability, maintaining nearly constant MSE values over the entire simulation duration. The magnitude of error is often smaller than that of FEM solutions computed with time step sizes 10 to 100 times finer. 
This superior performance shows the benefits of the discrete-time PINNs approach with a high-order IRK scheme ($q=50$) combined with the time-marching method in mitigating error propagation.

\begin{figure}
\centering
\includegraphics[width=1\textwidth]{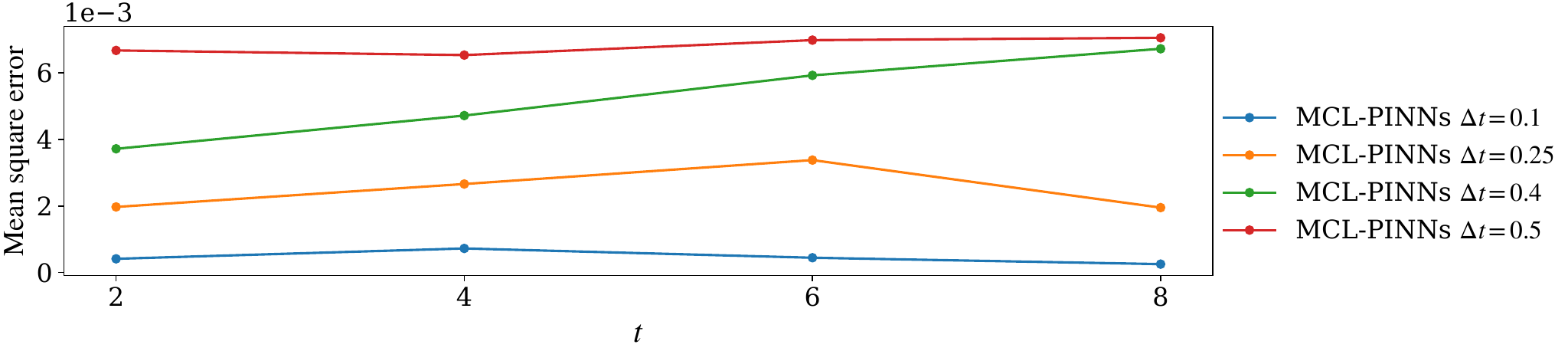}
\caption{Temporal evolution of mean squared error (MSE) for case $\theta_s = 50^\circ$ of the droplet coalescence problem under varying time step sizes $\Delta t = 0.1, 0.25, 0.4$, and $0.5$ in MCL-PINNs.}
\label{fig:mse_dt_50du_MCLPINNs}
\end{figure}

Figure~\ref{fig:mse_dt_50du_MCLPINNs} shows the sensitivity of the MCL-PINNs framework to the choice of time step size for the hydrophilic case ($\theta_s = 50^\circ$). Four different time step sizes are tested: $\Delta t = 0.1, 0.25, 0.4, 0.5$. The results demonstrate remarkable stability of the MCL-PINNs framework across all tested time step sizes. Even at $\Delta t = 0.5$, the mean squared error remains bounded and within a reasonable range throughout the simulation, while most  conventional numerical methods fail for this phase-field simulation at such a large time step size. For the selected time step size $\Delta t = 0.1$, the error converges to approximately $2.5 \times 10^{-4}$, suggesting this choice as an effective balance between computational cost and solution accuracy.
%------------------------------------------------------------------
%------------------------------------------------------------------
These findings indicate that MCL-PINNs can provide competitive predictive performance relative to conventional phase-field solvers, while showing potential advantages for long-time simulations.

%------------------------------------------------------------------
\subsection{Shear-induced droplet deformation} % case2 
In this subsection, we examine the deformation and migration of a droplet under shear flow, a classical problem involving moving contact lines. As noted in Remark \ref{remark1}, when the interface is driven by an external flow, an advection term is introduced into the Cahn--Hilliard equation. The governing equation for the order parameter $\phi$ becomes:
\begin{equation}
    \frac{\partial \phi}{\partial t} + \mathbf{u} \cdot \nabla \phi = M \Delta \mu,
\end{equation}
where $\mathbf{u} = (u_x, 0)$ is the prescribed velocity field. The computational domain is $\Omega = [0, 6] \times [0, 0.8]$ and $t \in [0, 8]$. We consider a Poiseuille flow profile for the velocity field $u_x(y) = C_0 \frac{4y(H-y)}{H^2}$, where $C_0$ represents the characteristic velocity magnitude. The initial condition describes a circular droplet centered at $(1, 0)$ with radius $r=0.4$:
\begin{equation}
    \phi^0(x, y) = \tanh\left( \frac{0.4 - \sqrt{(x - 1)^2 + y^2}}{2\varepsilon} \right).
\end{equation}
The following model parameters are used: $M=0.05$, $\varepsilon=0.02$, and  $\alpha=1000$. We examine two distinct wettability conditions characterized by the static contact angle $\theta_s$: a hydrophilic surface ($\theta_s = 50^\circ$) and a hydrophobic surface ($\theta_s = 150^\circ$). Dynamic boundary conditions \eqref{eq:relaxation_bc}-\eqref{eq:conservation_bc} are imposed on the bottom and top walls. Homogeneous Neumann boundary conditions $\nabla \phi \cdot \mathbf{n} = 0$ and $\nabla \mu \cdot \mathbf{n} = 0$ are applied on the left and right boundaries.

The proposed framework is configured with problem-specific hyperparameters. The variable-scaling transformation is applied with factors $N_x = 10$, $N_y = 10$, and $N_t = 1$, which balance the magnitudes of spatial and temporal derivatives. For the initial time interval, a uniform candidate grid of $375 \times 50$ nodes is constructed over the computational domain. In subsequent time intervals, this extraction procedure is applied to a refined candidate grid of $751 \times 101$ nodes. These interface-focused points are combined with additional interior samples distributed via LHS: 2,000 points for the hydrophilic case ($\theta_s = 50^\circ$) and 1,000 points for the hydrophobic case ($\theta_s = 150^\circ$). In addition, 158 uniformly distributed nodes are used for the horizontal boundaries and 21 nodes for the vertical boundaries. 
For the IRK scheme, the stage number is adjusted according to the numerical stiffness induced by the coupling of advection and wetting dynamics: $q = 50$ for $\theta_s = 50^\circ$ and $q = 100$ for $\theta_s = 150^\circ$. A fixed time step $\Delta t = 0.1$ is used. The number of training iterations is $N_{\text{iter}}=20,000$ for both contact angle configurations.

\begin{figure} 
    \centering
    % 第一行
    \begin{subfigure}[b]{0.48\textwidth}
        \centering
        \includegraphics[width=\linewidth]{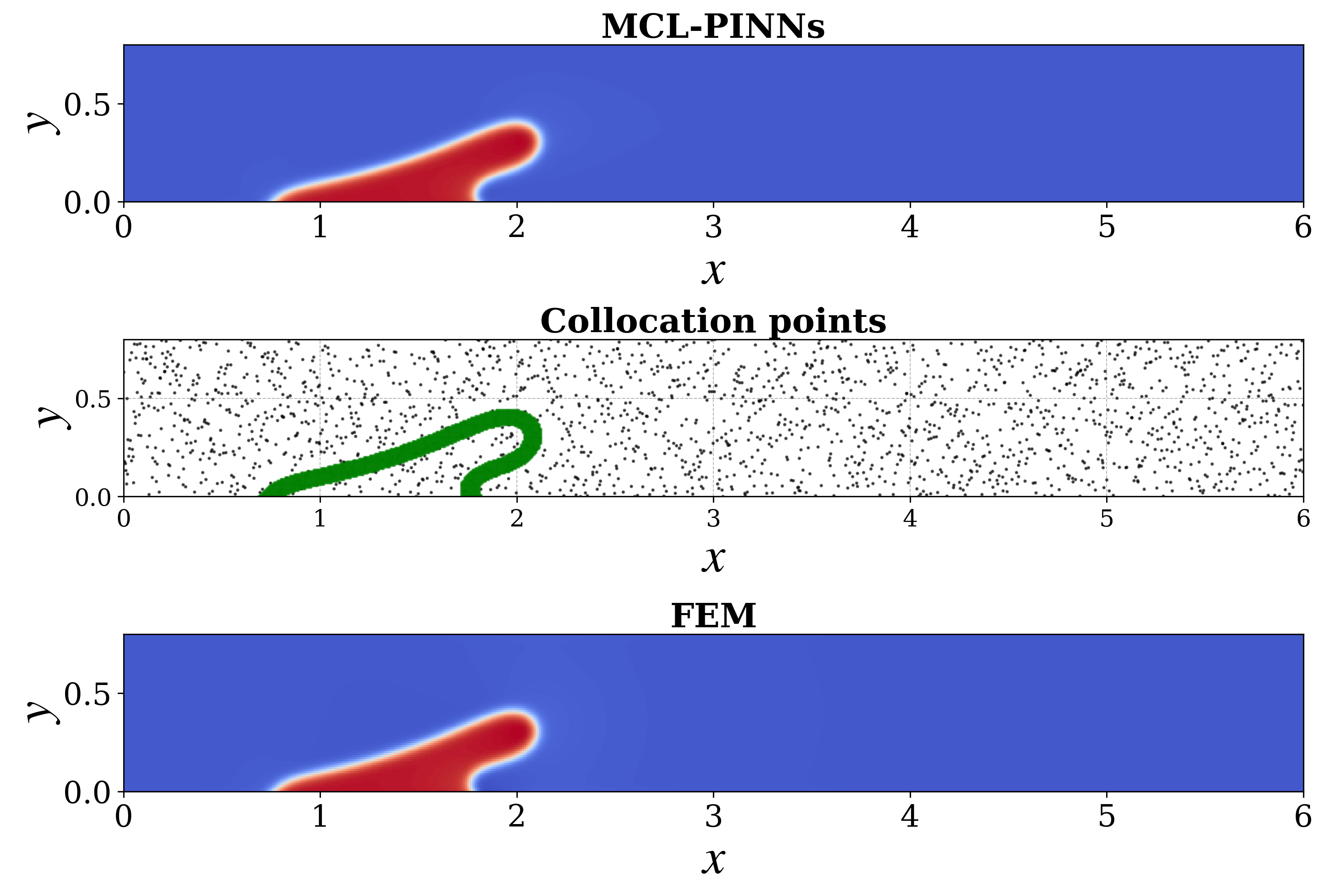}
        \caption{$\theta_s=50^\circ$, $T = 2$}
        \label{fig:theta50_T2}
    \end{subfigure}
    \hfill
    \begin{subfigure}[b]{0.48\textwidth}
        \centering
        \includegraphics[width=\linewidth]{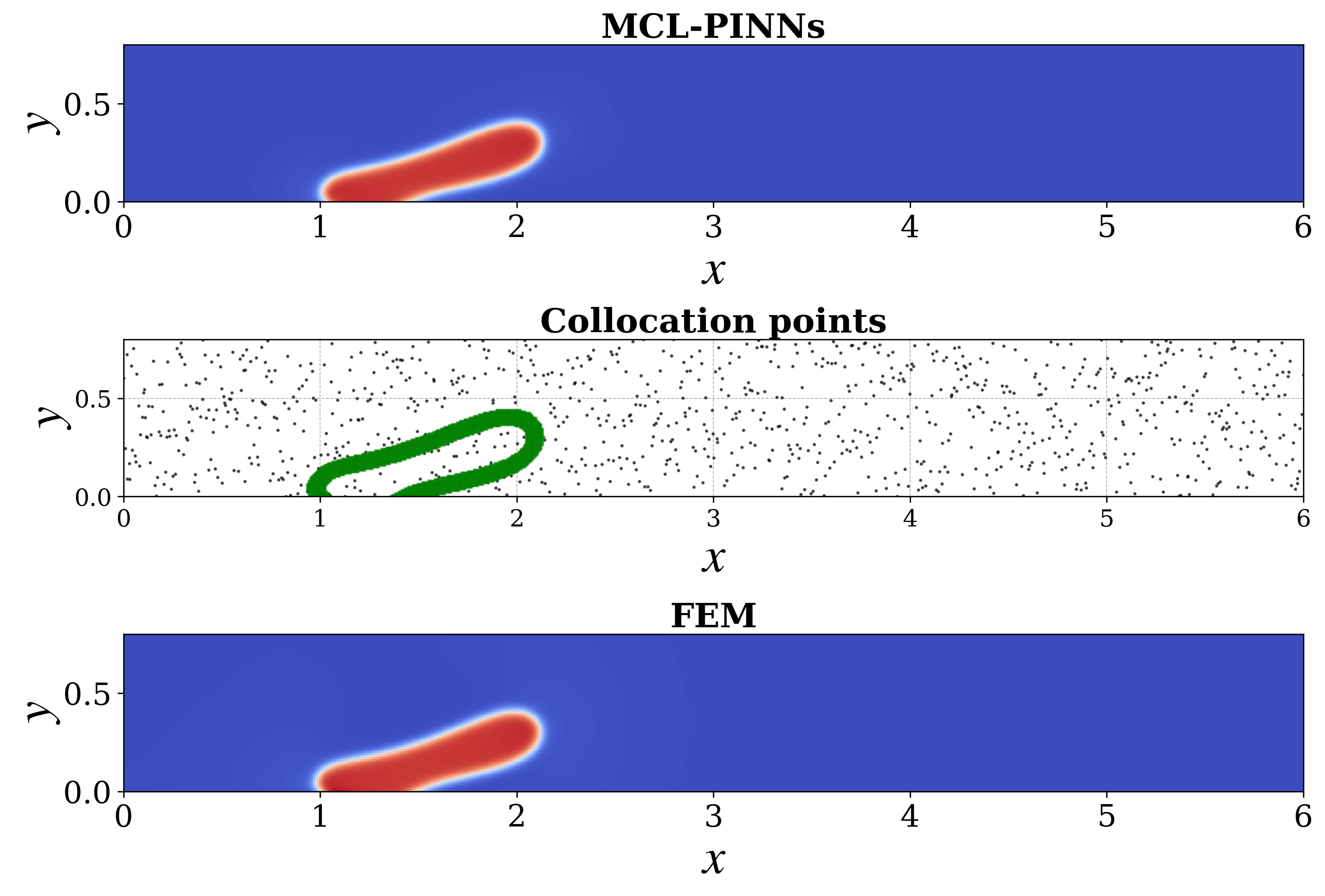}
        \caption{$\theta_s=150^\circ$, $T = 2$}
        \label{fig:theta150_T2}
    \end{subfigure}
      % \vspace{4em} 
    % 第二行
    \begin{subfigure}[b]{0.48\textwidth}
        \centering
        \includegraphics[width=\linewidth]{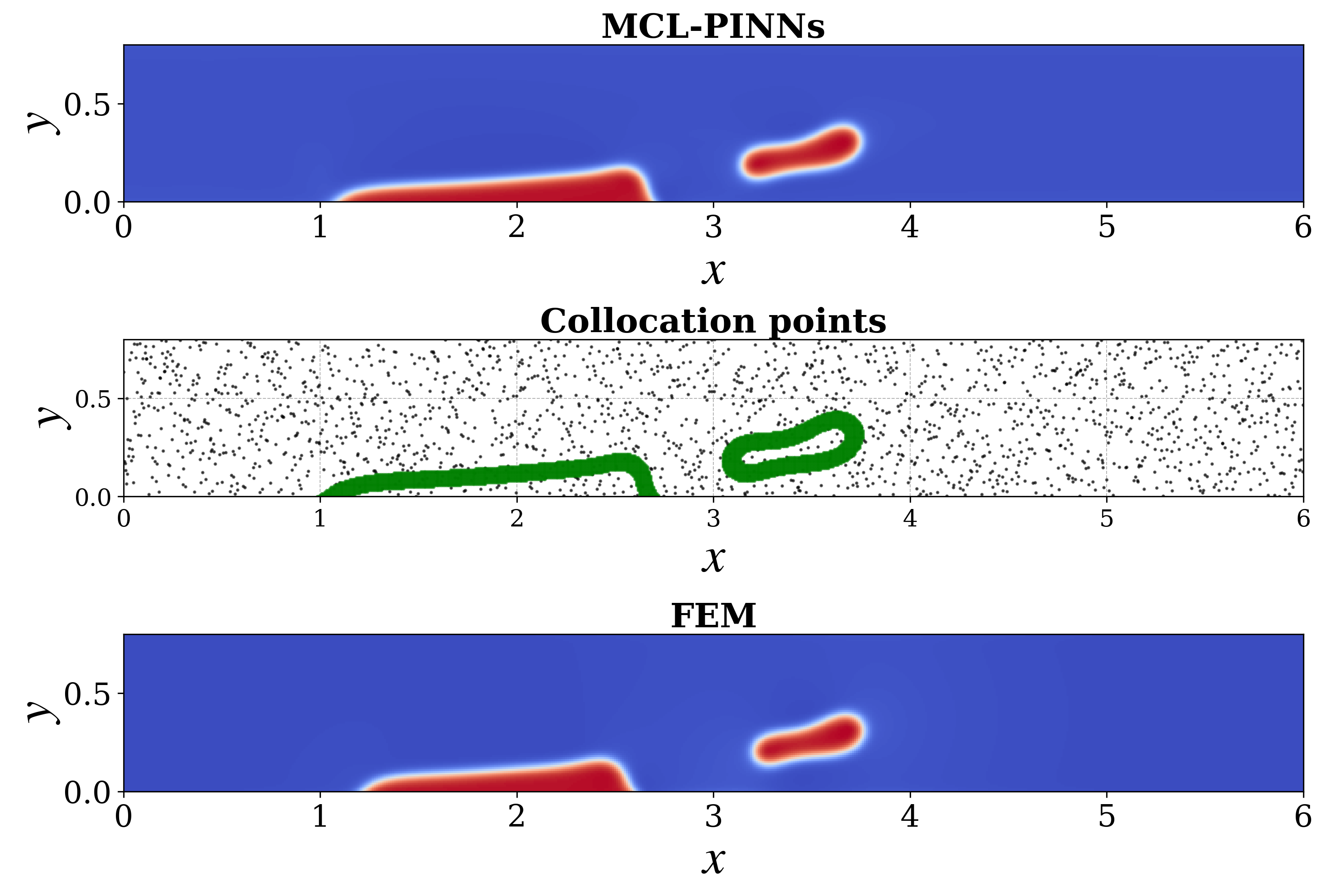}
        \caption{$\theta_s=50^\circ$, $T = 6$}
        \label{fig:theta50_T6}
    \end{subfigure}
    \hfill
    \begin{subfigure}[b]{0.48\textwidth}
        \centering
        \includegraphics[width=\linewidth]{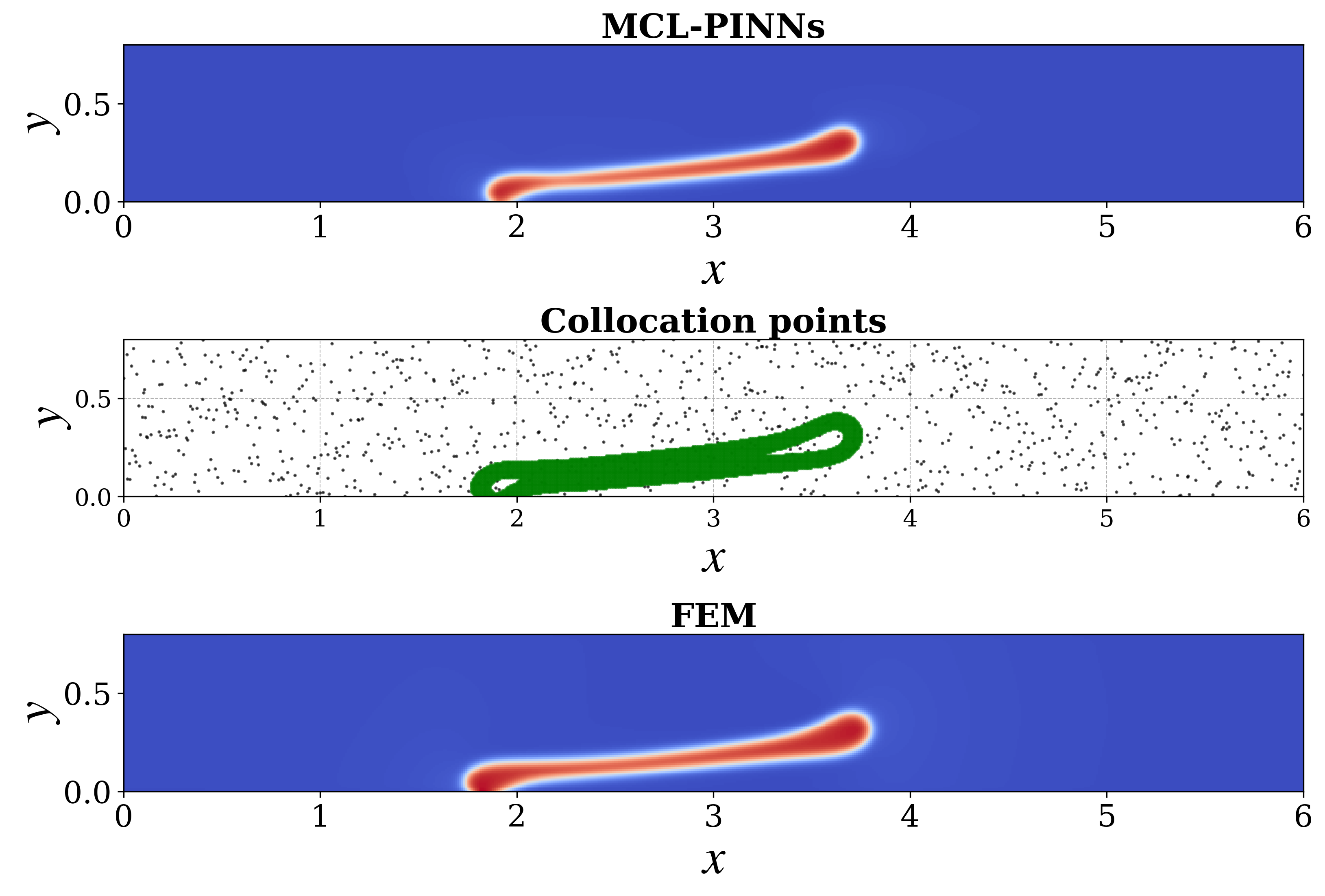}
        \caption{$\theta_s=150^\circ$, $T = 6$}
        \label{fig:theta150_T6}
    \end{subfigure}
    
    \caption{Comparison of phase field $\phi$ between MCL-PINNs and FEM solutions, along with the distribution of collocation points, for shear-induced droplet deformation at different static contact angles and time instances. The background flow is a scaled Poiseuille profile.}
    \label{fig:phase_field_comparison}
\end{figure}

First, we validate the effectiveness of the proposed MCL-PINNs by comparing its predictions with reference solutions obtained from the FEM. Figure~\ref{fig:phase_field_comparison} displays the phase-field variable $\phi$ at time instances $t=2$ and $t=6$ for a scaled flow velocity $C_0=0.5$. Under the imposed shear flow, the droplet interface undergoes pronounced deformation, accompanied by substantial contact line motion along the slip boundary.
For the hydrophilic case with $\theta_s=50^\circ$, as shown in Figures~\ref{fig:theta50_T2} and~\ref{fig:theta50_T6}, the strong affinity between the droplet and the solid substrate promotes spreading along the bottom boundary. As the shear flow drives the interface downstream, the upstream and downstream contact positions migrate at different rates, leading to a highly stretched interfacial profile. At the later time $T=6$, the upper portion of the droplet is further elongated and eventually separates from the near-wall liquid layer. In contrast, for the hydrophobic case with $\theta_s=150^\circ$, shown in Figures~\ref{fig:theta150_T2} and~\ref{fig:theta150_T6}, the weaker surface affinity limits the spreading of the droplet along the wall. The interface is primarily convected and elongated by the shear flow, while the contact line motion remains more localized compared with the hydrophilic case. In both wetting regimes, the MCL-PINNs predictions show good agreement with the FEM reference solutions, accurately capturing the interface position, curvature variation, and contact line evolution. The middle row of Figure~\ref{fig:phase_field_comparison} shows the corresponding spatial distribution of collocation points. The black dots denote the background samples generated by LHS, while the green dots represent the additional collocation points introduced by the adaptive sampling strategy described in Section~\ref{subsec:Adaptive_sampling}. As observed, the green points are concentrated along the diffuse interface and move consistently with the evolving interfacial region. This demonstrates that the adaptive sampling strategy effectively tracks the strongly deformed interface during the shear-induced wetting process. 

%------------------------------------------------------------------
%------------------------------------------------------------------

\begin{figure}
    \centering
    \includegraphics[width=0.8\linewidth]{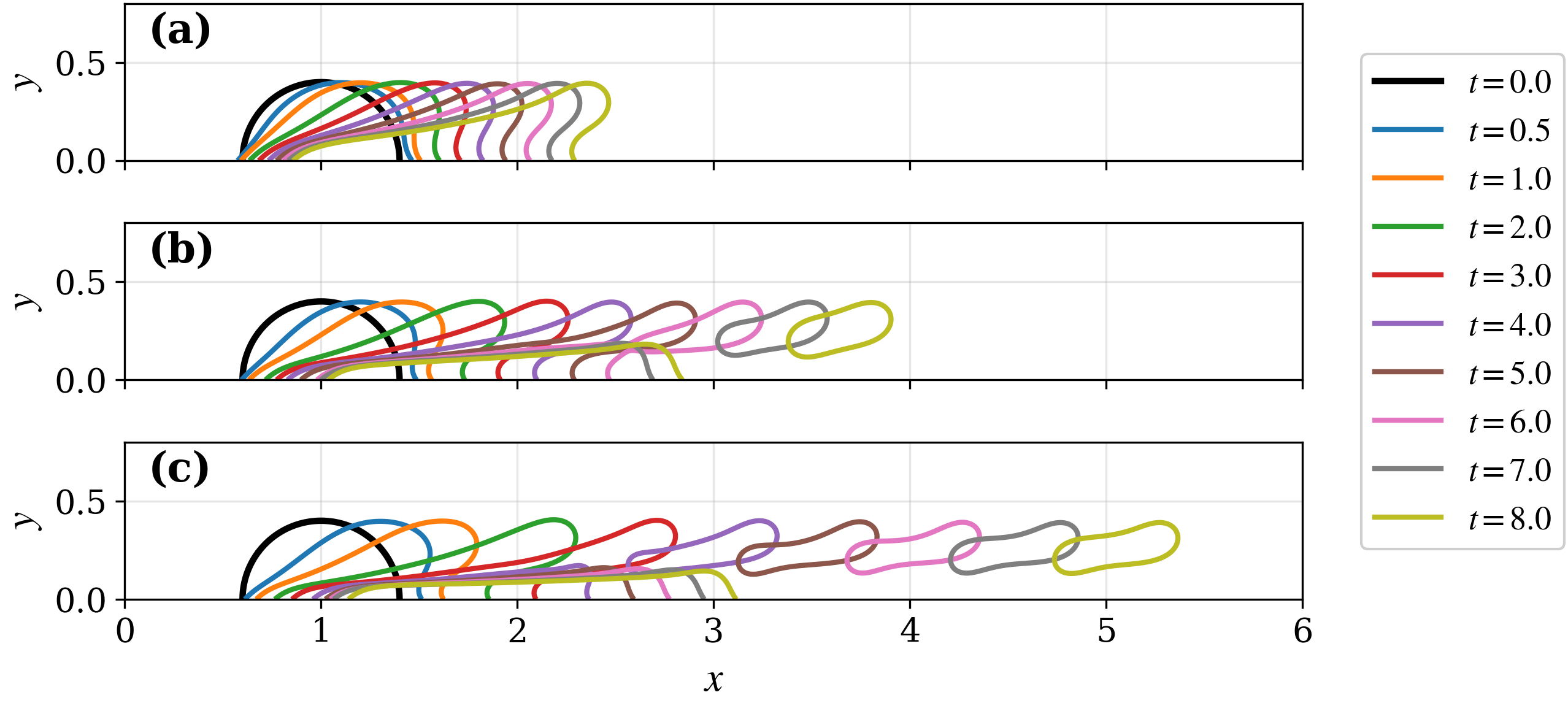}
    \caption{Evolution of the interface ($\phi=0$) for the hydrophilic case ($\theta_s=50^\circ$) at different flow intensities: (a) $C_0=0.2$, (b) $C_0=0.4$, and (c) $C_0=0.6$. The time range is from $t=0$ to $t=8$. Higher flow rates induce significant elongation and eventual pinch-off of the droplet.}
    \label{fig:case2_T=8_50du}
\end{figure}

\begin{figure}
    \centering
    \includegraphics[width=0.8\linewidth]{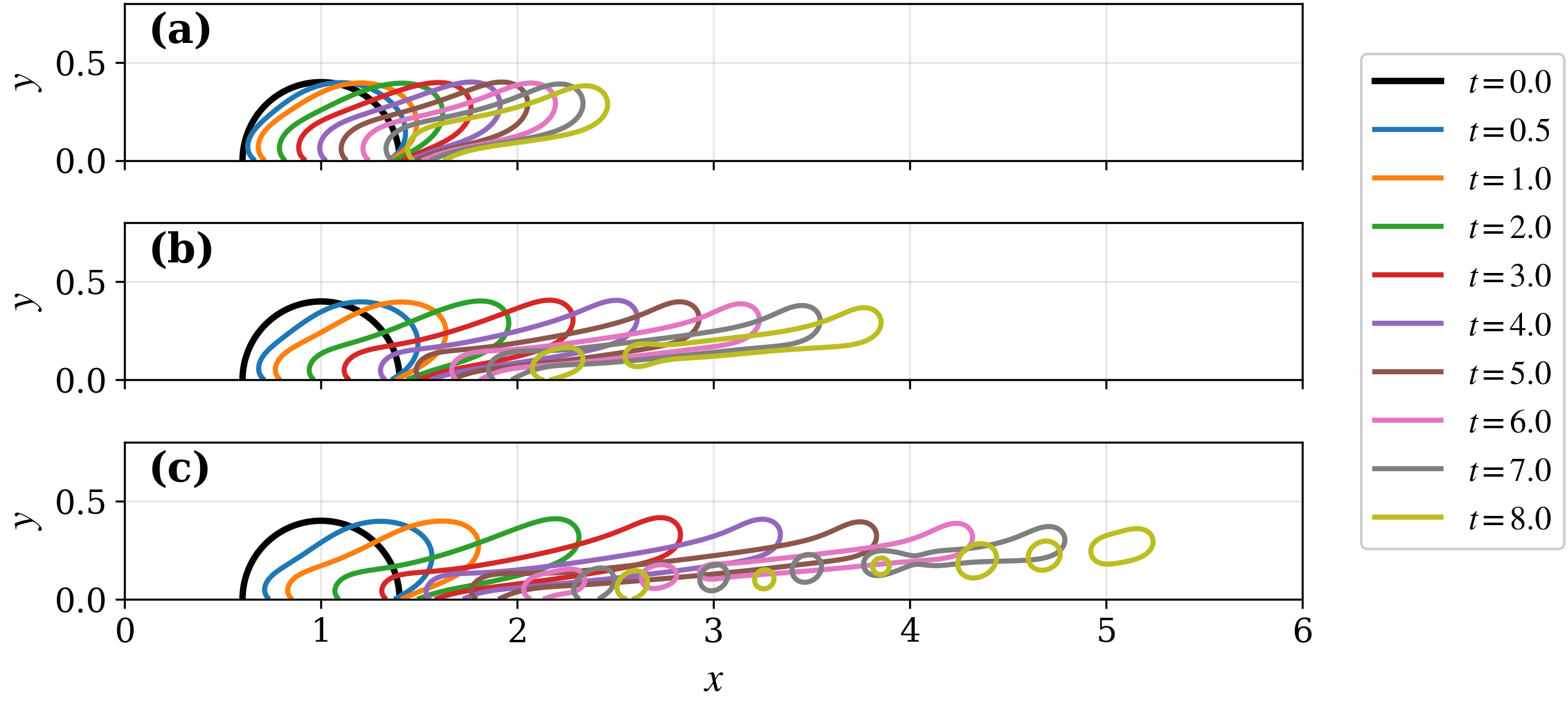}
    \caption{Evolution of the interface ($\phi=0$) for the hydrophobic case ($\theta_s=150^\circ$) at different flow intensities: (a) $C_0=0.2$, (b) $C_0=0.4$, and (c) $C_0=0.6$. The reduced wall adhesion leads to faster migration and easier breakup of the droplet compared to the hydrophilic case.}
    \label{fig:case2_T=8_150du}
\end{figure}

To further analyze the effects of flow velocity and wettability on droplet dynamics, we track the evolution of the interface ($\phi=0$) over the time interval $t\in[0,8]$. Figures~\ref{fig:case2_T=8_50du} and~\ref{fig:case2_T=8_150du} show the interface evolution under different flow velocities and wettability conditions. In each figure, subplots (a), (b), and (c) correspond to $C_0=0.2$, $C_0=0.4$, and $C_0=0.6$, respectively, while the static contact angle is set to $\theta_s=50^\circ$ in Figure~\ref{fig:case2_T=8_50du} and $\theta_s=150^\circ$ in Figure~\ref{fig:case2_T=8_150du}. As shown in these figures, increasing the convection velocity enhances the shear force, causing the upper portion of the droplet to move faster than the lower portion and thereby promoting droplet breakup. For cases with the same convection velocity but different contact angles, stronger wall hydrophobicity induces a more evident contraction of the lower portion of the droplet. This increases the volume of the droplet exposed to the high velocity region, leading to the formation of a larger detached droplet after pinch-off. These results demonstrate that the proposed MCL-PINNs can effectively capture the complex interplay between hydrodynamic forces, capillary forces, and wetting properties, accurately resolving topological changes such as droplet breakup.
%------------------------------------------------------------------
%------------------------------------------------------------------
%------------------------------------------------------------------

%------------------------------------------------------------------
\subsection{Dynamic wetting through a heterogeneous channel} \label{sec:het_wetting} % case3 
\begin{figure}
    \centering
    \includegraphics[width=0.7\textwidth, height=0.25\textheight, keepaspectratio]{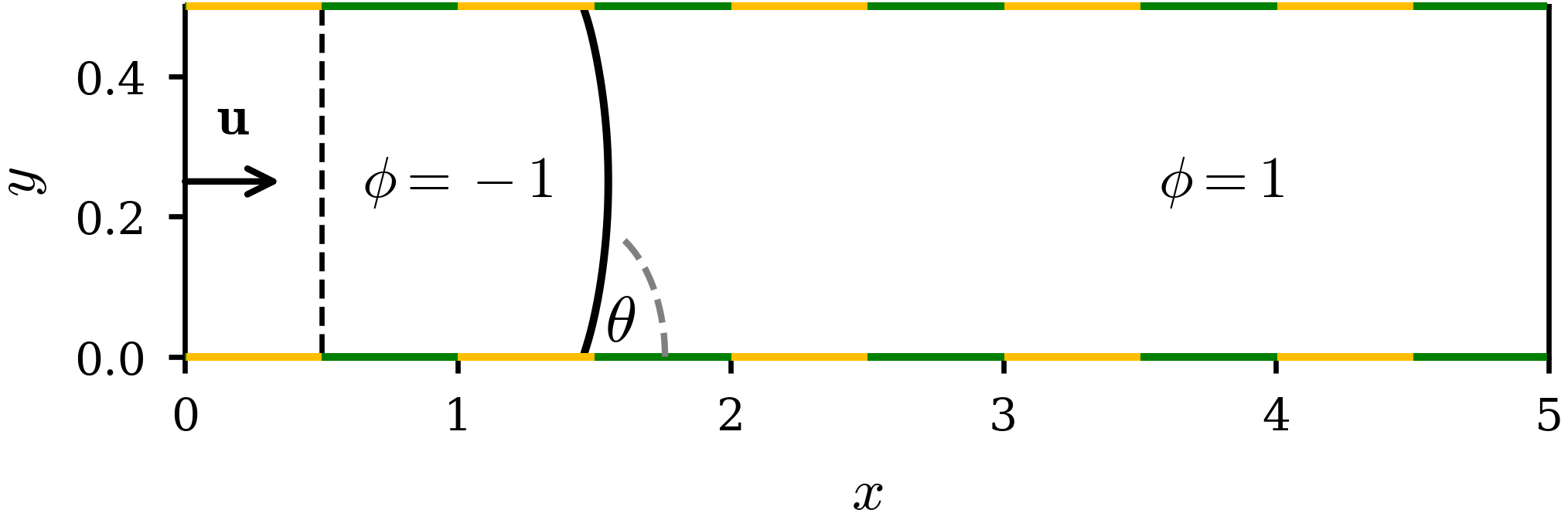}
    \caption{An illustration of the computational domain and boundary configuration for the heterogeneous channel problem. The top and bottom boundaries are periodically patterned with two wetting properties, where yellow and green segments denote regions with static contact angles $\theta_{s_1}=60^\circ$ and $\theta_{s_2}=120^\circ$, respectively. The black dashed line denotes the initial position of the interface.}
    \label{fig:het_schematic}
\end{figure}

In this subsection, we investigate the dynamic wetting behavior of a binary fluid interface driven over a chemically patterned channel. 
An illustration of the problem setup is shown in Figure~\ref{fig:het_schematic}. The computational domain is $\Omega = [0, 5] \times [0, 0.5]$ and $t \in [0, 9]$. The top and bottom boundaries are chemically patterned in a periodic manner, alternating between two wetting properties. Specifically, the dynamic boundary conditions~\eqref{eq:relaxation_bc}-\eqref{eq:conservation_bc} are imposed with a spatially varying static contact angle $\theta_s(x)$ with period $1$:
\begin{equation} \label{eq:theta_s_periodic}
\theta_s(x)=\Theta(\xi), 
\qquad 
\xi=x-\lfloor x \rfloor \in [0,1),
\end{equation}
where
\begin{equation}
\Theta(\xi)=
\begin{cases}
\theta_{s_2} - \dfrac{(\xi+\delta)(\theta_{s_2}-\theta_{s_1})}{2\delta},
& 0 \leq \xi < \delta, \\
\theta_{s_1},
& \delta \leq \xi < \dfrac{1}{2}-\delta, \\
\theta_{s_1} + \dfrac{(\xi-\frac{1}{2}+\delta)(\theta_{s_2}-\theta_{s_1})}{2\delta},
& \dfrac{1}{2}-\delta \leq \xi < \dfrac{1}{2}+\delta, \\
\theta_{s_2},
& \dfrac{1}{2}+\delta \leq \xi < 1-\delta, \\
\theta_{s_2} - \dfrac{(\xi-1+\delta)(\theta_{s_2}-\theta_{s_1})}{2\delta},
& 1-\delta \leq \xi < 1.
\end{cases}
\end{equation}
Here, $\theta_{s_1}=60^\circ$, $\theta_{s_2}=120^\circ$, and the transition half-width is $\delta=\varepsilon/2=0.0125$.
Homogeneous Neumann boundary conditions $\nabla \phi \cdot \mathbf{n} = 0$ and $\nabla \mu \cdot \mathbf{n} = 0$ are applied on the left and right boundaries. A uniform velocity field $\mathbf{u} = (u_x, 0)$ is prescribed as a driving force, where $u_x = 0.5$ for advancing motion and $u_x = -0.5$ for receding motion. The initial condition is given as $\phi^0(x, y) = \tanh\bigl((x - x_c)/(2\varepsilon)\bigr)$, where $x_c = 0.5$ for the advancing case and $x_c = 4.5$ for the receding case. The following model parameters are used: $M=1$, $\varepsilon=0.025$, and $\alpha=100$. 
%------------------------------------------------------------------
%------------------------------------------------------------------
For the proposed MCL-PINNs framework, the variable-scaling transformation is applied with scaling factors $N_x=10$, $N_y=10$, and $N_t=1$. A uniform $751\times 101$ grid is used as the candidate pool for adaptive sampling. These points are further supplemented by $2,000$ interior collocation points generated using LHS to provide baseline coverage of the computational domain. 
For the boundary residuals, $158$ uniformly distributed nodes are placed on the horizontal boundaries, while $21$ nodes are used on the vertical boundaries. In both the advancing and receding configurations, we use a fixed time step $\Delta t=0.1$ and an IRK scheme with $q=50$ stages. The number of training iterations is set to $N_{\text{iter}}=10,000$ for both cases.

\begin{figure}
    \centering
    \includegraphics[width=0.7\textwidth, height=0.25\textheight, keepaspectratio]{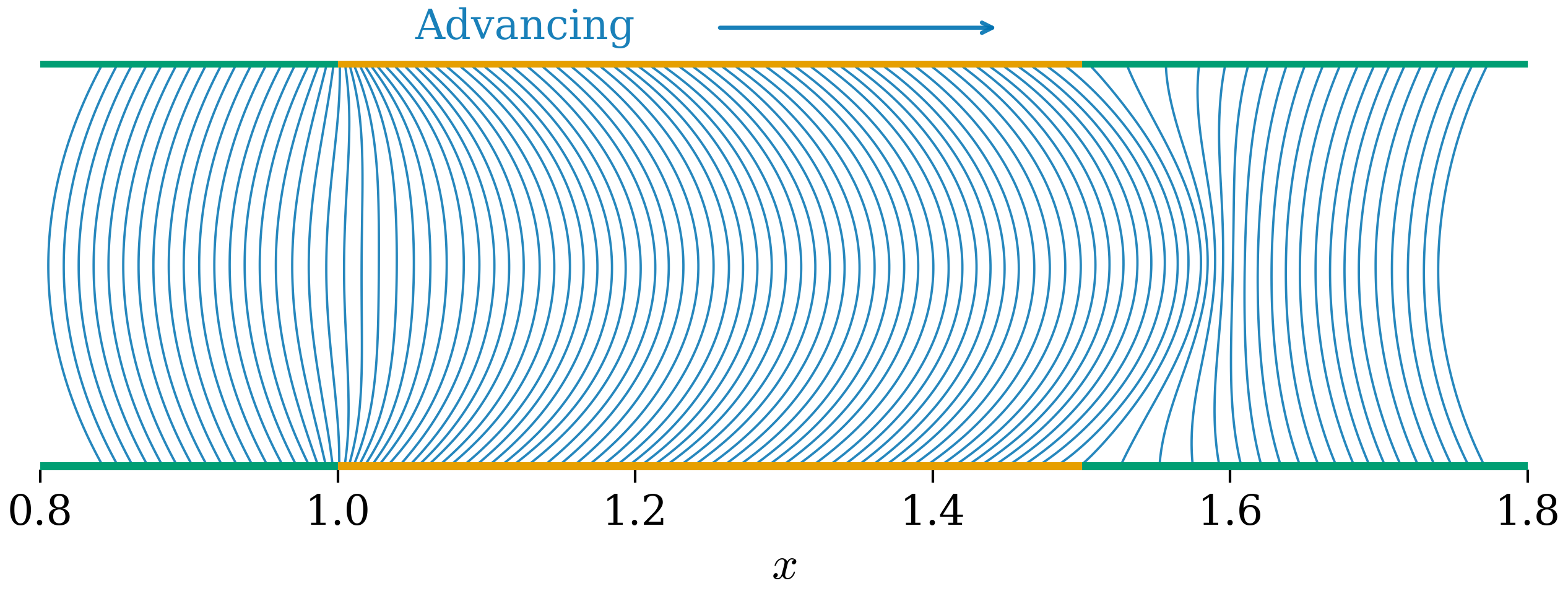}
    \caption{
Advancing phase interface in a channel with heterogeneous wetting patterns. The blue curves represent the phase interfaces at successive time instants, with a time interval of $\Delta t=0.02$ between adjacent curves. The alternating yellow and green wall segments denote hydrophilic regions with $\theta_{s_1}=60^\circ$ and hydrophobic regions with $\theta_{s_2}=120^\circ$, respectively. The interface is driven from left to right by a uniform velocity $u_x=0.5$.
        }
    \label{fig:advancing_heterogeneous}
\end{figure}

%------------------------------------------------------------------
Figure~\ref{fig:advancing_heterogeneous} shows the interfacial motion through a heterogeneous channel under the imposed driving force. The blue curves represent the phase interfaces at successive time instants, with a time interval of $\Delta t=0.02$ between adjacent profiles. As the contact points pass across the junctions between regions with different wetting properties, the interface undergoes pronounced deformation and exhibits characteristic \textit{stick--slip} motion.
When the phase $\phi=-1$ advances from the green region to the yellow region, it moves toward a less favorable wetting region and therefore experiences a larger surface-energy barrier. Consequently, the contact line is temporarily pinned near the top and bottom chemical junctions, leading to slow displacement of the contact points. During this stage, the interface is stretched and capillary energy gradually accumulates. Once the hydrodynamic driving force overcomes the local pinning threshold, depinning occurs and the contact line rapidly slips forward.
In contrast, when the phase $\phi=-1$ moves from the yellow region to the green region, the motion is directed toward a more favorable wetting region. The resistance at the chemical junction is therefore reduced, allowing the interface to pass through the junction more rapidly and resulting in a slip-dominated response. 

%------------------------------------------------------------------

\begin{figure}
    \centering
    \includegraphics[width=\linewidth]{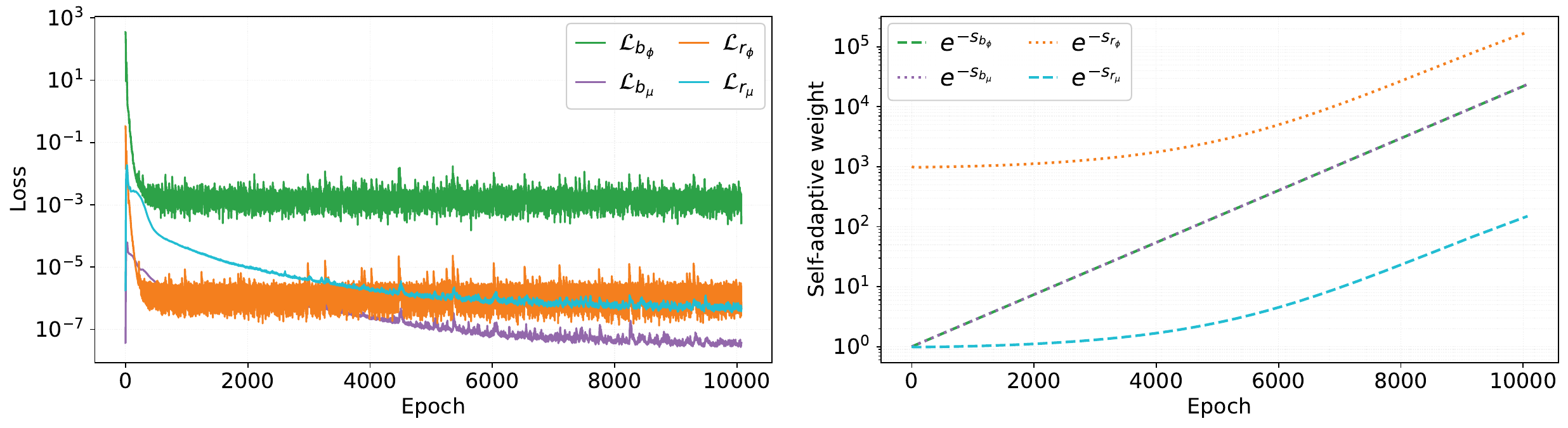}
    \caption{
    (Left) Logarithmic evolution of the boundary and residual loss components. 
    (Right) Corresponding evolution of the adaptive weights, which dynamically adjust to balance the magnitude of each loss term.}
    \label{fig:training_dynamics_case3_adv}
\end{figure}

To understand how the objective is optimized during training, we examine the evolution of the individual loss components and the corresponding learned adaptive weights over an intermediate time interval, $t\in[0.3,0.4]$, as shown in Figure~\ref{fig:training_dynamics_case3_adv}. The left panel shows that all residual terms decrease rapidly in the early stage and continue to decay over subsequent epochs, although their magnitudes remain different. In particular, the boundary residual $\mathcal{L}_{b_\phi}$ remains relatively larger, reflecting the stronger stiffness and nonlinear coupling associated with the dynamic boundary condition for $\phi$. In contrast, $\mathcal{L}_{r_\phi}$, $\mathcal{L}_{r_\mu}$, and $\mathcal{L}_{b_\mu}$ are reduced to much smaller levels as training proceeds. The right panel presents the evolution of the learned self-adaptive weights $e^{-s_m}$ for the individual loss components. These weights are optimized during training instead of being manually prescribed. The different growth rates indicate that the adaptive weighting strategy assigns different levels of emphasis to the residual terms according to their optimization behavior. In particular, $e^{-s_{r_\phi}}$ exhibits the largest value, showing that the phase-field residual receives a progressively stronger weight in the objective. The boundary-condition weights $e^{-s_{b_\phi}}$ and $e^{-s_{b_\mu}}$ also increase substantially, which strengthens the enforcement of the dynamic boundary conditions. In comparison, $e^{-s_{r_\mu}}$ is relatively small, indicating that the chemical-potential residual is more readily balanced during training. 

\begin{figure}
    \centering
    \includegraphics[width=0.8\textwidth, height=0.65\textheight, keepaspectratio]{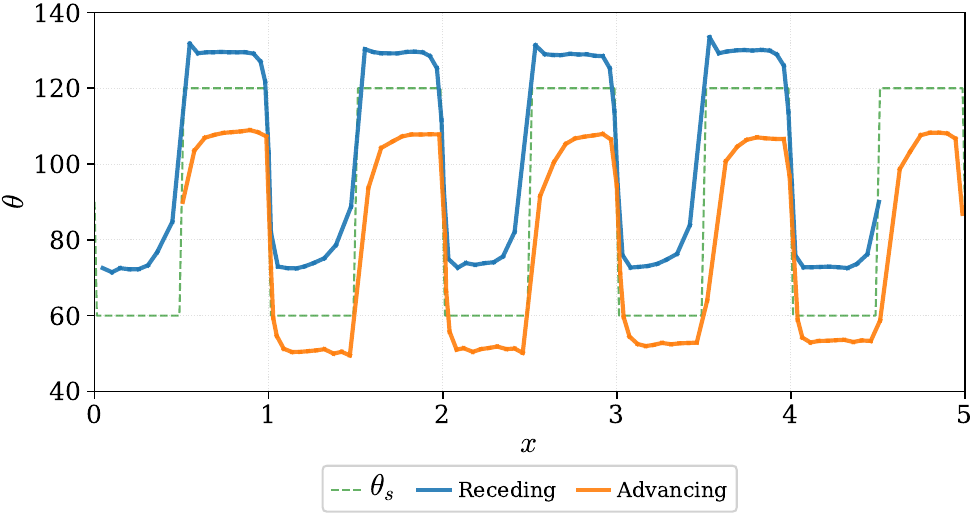}
    \caption{Evolution of the dynamic contact angle $\theta$ along the channel walls as a function of the contact line position $x$. The solid curves represent the computed contact angles during advancing and receding motions, while the dashed line indicates the prescribed static contact angle $\theta_s(x)$. The difference between the advancing and receding curves indicates a clear contact angle hysteresis induced by the heterogeneous walls. }
    \label{fig:het_angle}
\end{figure}

\begin{figure}
    \centering
    \includegraphics[width=\linewidth]{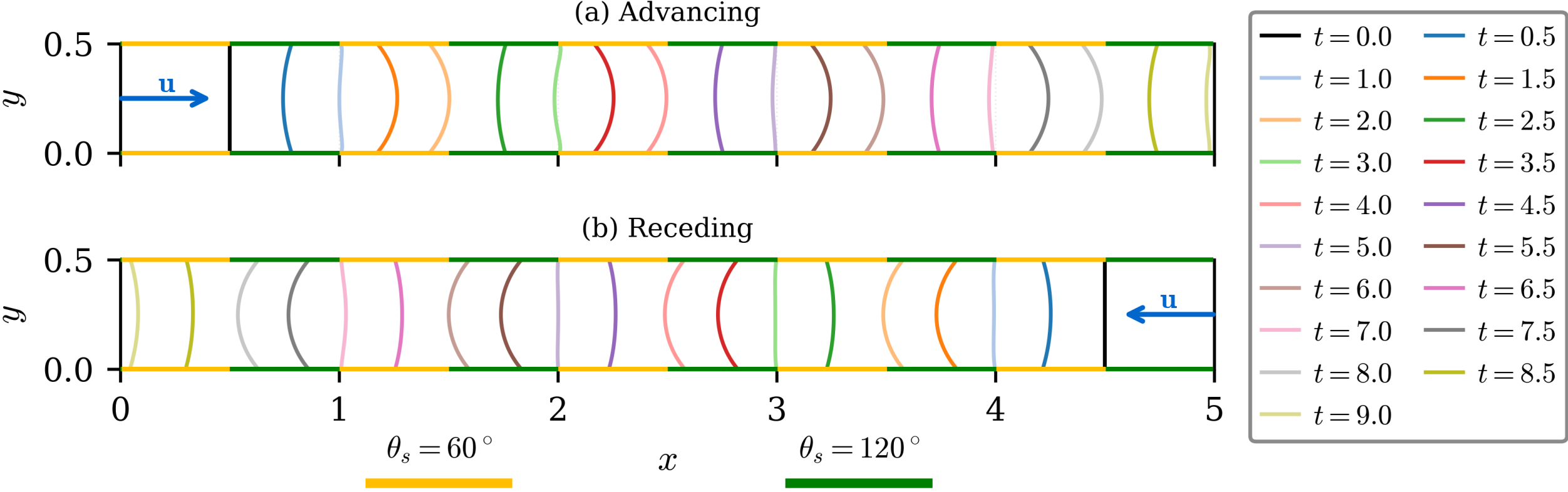}
    \caption{Comparison of the interface ($\phi=0$) trajectories between advancing and receding cases. 
    }
    \label{fig:het_phi_zero}
\end{figure}

The quantitative wetting response is further analyzed in Figure~\ref{fig:het_angle}, where the dynamic contact angle $\theta$ is plotted as a function of the contact position $x$. The corresponding interfacial evolution during the advancing and receding processes is shown in Figure~\ref{fig:het_phi_zero}. In Figure~\ref{fig:het_angle}, the orange curve corresponds to the advancing case in Figure~\ref{fig:het_phi_zero}(a), whereas the blue curve corresponds to the receding case in Figure~\ref{fig:het_phi_zero}(b). It should be noted that the computed dynamic contact angle does not coincide with the prescribed static contact angle profile $\theta_s(x)$. The static contact angle represents the local equilibrium wettability of the channel walls, whereas the dynamic contact angle is an apparent non-equilibrium quantity determined by the instantaneous interface morphology, contact-line motion, and relaxation of the phase field near the wall. Therefore, when the contact line moves across heterogeneous wettability patterns, the interface cannot instantaneously adjust to the local value of $\theta_s(x)$, leading to a delayed and path-dependent response. Moreover, a clear contact angle hysteresis induced by the heterogeneous walls can be observed. In particular, several sharp variations in the contact angle occur near the wettability transition points, for example at $x=0.5$, $1.5$, $2.5$, and $3.5$ in the receding case. When the interface moves from the green boundary region to the yellow boundary region, the contact line encounters stronger surface-energy resistance, since the phase $\phi=-1$ being displaced tends to remain in the energetically favorable region. As a result, the angle between the interface and the solid boundary increases abruptly. Subsequently, driven by interfacial-energy relaxation, the contact angle gradually decreases toward a smaller value. A similar mechanism can be observed in the advancing case, where the response follows a different branch due to the opposite direction of contact-line motion. These results demonstrate that the heterogeneous wettability pattern gives rise to a pronounced advancing-receding asymmetry and a clear hysteresis loop in the dynamic contact angle response. The close agreement between the dynamic contact angle response and the corresponding interface evolution further demonstrates that the proposed MCL-PINNs can effectively capture the hysteretic wetting behavior induced by boundary heterogeneity.

%------------------------------------------------------------------
%------------------------------------------------------------------

%-------------------------------------------------------------------------

\section{Concluding remarks} \label{sec:conclusion}
In this study, we propose MCL-PINNs, a neural solver for moving contact line problems with dynamic boundary conditions. By embedding the Cahn--Hilliard equation and the associated boundary constraints into a physics-informed neural network framework, the proposed method addresses several numerical challenges arising from long-time error accumulation, sharp interfacial profiles, localized contact line dynamics, and complex contact angle evolution.
Building upon the discrete-time formulation, the proposed framework integrates several key techniques, including a time-marching strategy with multiple networks, a relaxed distribution constraint on the neural network outputs, a variable-scaling method, adaptive loss weighting, adaptive collocation sampling, and, when applicable, symmetry preservation through neural network inputs. These components work together to improve training stability, enhance the balance among competing loss terms, and accurately capture complex interfacial evolution in moving contact line dynamics.

Numerical validation is conducted on three examples, including droplet coalescence, shear-induced droplet deformation, and dynamic wetting in a heterogeneous channel. The results demonstrate the effectiveness and robustness of the proposed method. The computed solutions show good agreement with reference solutions obtained by finite element methods, indicating that MCL-PINNs can accurately capture the governing phase-field dynamics and the associated moving contact line motion. To the best of our knowledge, this work represents the first PINN-type approaches for solving the Cahn--Hilliard equation with dynamic boundary conditions, providing a step toward connecting mesh-free deep learning solvers with thermodynamically consistent phase-field modeling of moving contact line problems.
One direction for future research is to extend the present framework to the coupled Cahn--Hilliard--Navier--Stokes (CHNS) system for moving contact line problems with dynamic boundary conditions. 

% %Bibliography
\bibliographystyle{unsrt}  
\bibliography{refs}

\end{document}